\documentclass[12pt]{article}
\usepackage{amssymb,amsmath,amsfonts,bbm,pifont,upgreek,bbold,accents,color}  
\usepackage{hyperref}
 
\font\teneufm=eufm10
\font\seveneufm=eufm7
\font\fiveeufm=eufm5
\newfam\eufmfam
\textfont\eufmfam=\teneufm
\scriptfont\eufmfam=\seveneufm
\scriptscriptfont\eufmfam=\fiveeufm

\newcommand\beq[1]{ \begin{equation}\label{#1} }
\newcommand{\eeq}{ \end{equation} }

\newcommand{\beqa}{ \begin{align*} }
\newcommand{\eeqa}{ \end{align*} }

\newcommand{\beqano}{ \begin{eqnarray*} }
\newcommand{\eeqano}{ \end{eqnarray*} }
\newtheorem{theorem}{Theorem}
\newtheorem{definition}[theorem]{Definition}
\newtheorem{proposition}{Proposition}
\newtheorem{lemma}[theorem]{Lemma}
\newtheorem{sublemma}{Sublemma}
\newtheorem{remark}[theorem]{Remark}
\newtheorem{notationalremark}{Notations}
\newtheorem{corollary}[theorem]{Corollary}
\newtheorem{assumption}{Assumption}
\newtheorem{claim}{Claim}
\newcommand{\zerarcounters}{\setcounter{equation}{0}\setcounter{theorem}{0}}

\numberwithin{theorem}{section}

\newtheorem{tools}{$\negsp\negsp$}[subsection]

\newcommand{\T}{\mathcal{T}}
\newcommand{\R}{\mathcal{R}}

\newcommand\thm[1]{ \begin{theorem}\label{#1}}
\newcommand\thmtwo[2]{ \begin{theorem}[#1]\label{#2}}
\newcommand\ethm{ \end{theorem} }
\newcommand\dfn[1]{ \begin{definition}\label{#1} \rm}
\newcommand\dfntwo[2]{ \begin{definition}[#1]\label{#2} \rm}
\newcommand\edfn{ \end{definition} }
\newcommand\pro[1]{ \begin{proposition}\label{#1}}
\newcommand\protwo[2]{ \begin{proposition}[#1]\label{#2}}
\newcommand\epro{ \end{proposition} }
\newcommand\lem[1]{ \begin{lemma}\label{#1}}
\newcommand\lemtwo[2]{ \begin{lemma}[#1]\label{#2}}
\newcommand\elem{ \end{lemma} }
\newcommand\sublem[1]{ \begin{sublemma}\label{#1}}
\newcommand\sublemtwo[2]{ \begin{sublemma}[#1]\label{#2}}
\newcommand\esublem{ \end{sublemma} }
\newcommand\rem[1]{ \begin{remark}\label{#1} \rm}
\newcommand\erem{ \end{remark} }
\newcommand\notrem[1]{ \begin{notationalremark}\label{#1} \rm}
\newcommand\enotrem{ \end{notationalremark} }
\newcommand\cor[1]{ \begin{corollary}\label{#1}}
\newcommand\cortwo[2]{ \begin{corollary}[#1]\label{#2}}
\newcommand\ecor{ \end{corollary} }
\newcommand\asmp[1]{ \begin{assumption}\label{#1}}
\newcommand\asmptwo[2]{ \begin{assumption}[#1]\label{#2}}
\newcommand\easmp{ \end{assumption} }
\newcommand\clm[1]{ \begin{claim}\label{#1}}
\newcommand\eclm{ \end{claim} }
\newcommand{\proof}{\par\medskip\noindent{\bf Proof\ }}
\expandafter\chardef\csname pre amssym.def
at\endcsname=\the\catcode`\@
\catcode`\@=11
\def\undefine#1{\let#1\undefined}
\def\newsymbol#1#2#3#4#5{\let\next@\relax
 \ifnum#2=\@ne\let\next@\msafam@\else
 \ifnum#2=\tw@\let\next@\msbfam@\fi\fi
 \mathchardef#1="#3\next@#4#5}
\def\mathhexbox@#1#2#3{\relax
 \ifmmode\mathpalette{}{\m@th\mathchar"#1#2#3}%
 \else\leavevmode\hbox{$\m@th\mathchar"#1#2#3$}\fi}
\def\hexnumber@#1{\ifcase#1 0\or 1\or 2\or 3\or 4\or 5\or 6\or 7\or
8\or
 9\or A\or B\or C\or D\or E\or F\fi}
\ifcase\@ptsize
 \font\tenmsb=msbm10
 \font\sevenmsb=msbm7
 \font\fivemsb=msbm5
\or
 \font\tenmsb=msbm10 scaled \magstephalf
 \font\sevenmsb=msbm7 scaled \magstephalf
 \font\fivemsb=msbm5  scaled \magstephalf
\or
 \font\tenmsb=msbm10 scaled \magstep1
 \font\sevenmsb=msbm7 scaled \magstep1
 \font\fivemsb=msbm5 scaled \magstep1
\fi
\newfam\msbfam
\textfont\msbfam=\tenmsb
\scriptfont\msbfam=\sevenmsb
\scriptscriptfont\msbfam=\fivemsb
\edef\msbfam@{\hexnumber@\msbfam}
\def\Bbb#1{\fam\msbfam\relax#1}
\def\widehat#1{\setboxz@h{$\m@th#1$}%
 \ifdim\wdz@>\tw@ em\mathaccent"0\msbfam@5B{#1}%
 \else\mathaccent"0362{#1}\fi}
\def\widetilde#1{\setboxz@h{$\m@th#1$}%
 \ifdim\wdz@>\tw@ em\mathaccent"0\msbfam@5D{#1}%
 \else\mathaccent"0365{#1}\fi}

\def\RIfM@{\relax\ifmmode}
\def\nonmatherr@#1{\errmessage{\string#1\space allowed only in math mode}}
\def\Bbb{\RIfM@\expandafter\Bbb@\else
 \expandafter\nonmatherr@\expandafter\Bbb\fi}
\def\Bbb@#1{{\Bbb@@{#1}}}
\def\Bbb@@#1{\fam\msbfam\relax#1}
\def\setboxz@h{\setbox\z@\hbox}
\def\wdz@{\wd\z@}
\catcode`\@=\csname pre amssym.def at\endcsname
\newcommand{\giu}{{\medskip\noindent}}

\newcommand{\nl}{{\smallskip\noindent}}

 \newcommand{\Fullbox}{{\rule{2.0mm}{2.0mm}}} 
 \newcommand{\qed}{\hfill\Fullbox\vspace{0.2cm}} 
\newcommand{\negsp}{\hspace{-.09truecm}}  %%% equivalente a \!

\newcommand{\dst}{\displaystyle}

\newcommand\su[1]{ \frac{1}{ {#1}} }

\renewcommand{\a }{ {\alpha}   }
\renewcommand{\b}{ {\beta}   }

\newcommand{\e }{ {\epsilon}   }

\newcommand{\x }{ {\xi}   }

\renewcommand{\t}{ {\tau}   }

\newcommand\by{{ \bar y}}

\newcommand{\gotd}{{\mathfrak d}}

\renewcommand\subset{\subseteq}

\newcommand{\Diff}{{\rm Diff}}

\begin{document}

\title{\bf On the linearization of analytic diffeomorphisms of the torus}

\author
{\bf Fernando Argentieri$^{1}$, Livia Corsi$^{2}$
\vspace{2mm}
\\ \small
$^{1}$ Institut f\"ur Mathematik, Universit\"at Z\"urich
Winterthurerstrasse 190, CH-8057 Z\"urich, CH
\\ \small 
$^{2}$ Dipartimento di Matematica, Universit\`a di Roma Tre, Roma,
I-00146, Italy
\\ \small 
E-mail: fernando.argentieri@math.uzh.ch, livia.corsi@uniroma3.it}

\date{}

\maketitle

\begin{abstract}
We provide an arithmetic condition weaker then the Bryuno condition for which it is possible to apply
a KAM scheme in dimension greater then one. The KAM scheme will be provided in the setting of 
linearization of analytic diffeomorphisms of the torus that are close to a rotation. 

\smallskip


\smallskip

%\noindent{\bf MSC classification}: 		37J40, 	70H08, 37C55

\end{abstract}

\section{Introduction and main results}

The study of existence and stability of invariant tori for quasi-integrable hamiltonian systems is one of the central problems
in dynamics, and origins from the three body problem. Due to the seminal work of Poincar\'e, it
was widely believed that in general invariant tori would not survive 
a small perturbation, beacuse of the obstruction due to accumulation of small divisors. 

Nevertheless, independently, Siegel in \cite{sie} and Kolmogorov in \cite{kol} were able to bypass the presence of small 
divisors\footnote{In particular Siegel considered the problem of linearization of analytic functions near elliptic fixed points, whereas Kolmogorov considered
 the problem of existence and stability of invariant tori.} by imposing some arithmetic condition on the frequency:
specifically they both imposed a diophantine condition. Then the question arises wheter imposing an arithmetic condition on the frequency is
necessary to overcome the presence of small divisors, or it is simply a technical tool\footnote{The fact that small divisors are really an 
obstruction for linearization in the smooth case can be seen already in the Anosov-Katok construction first introduced in \cite{AK}. 
More recently in \cite{FK} Fayad and Katok implemented for the first time the Anosov-Katok construction in the analytic case.}.

The  non-optimality of the diophantine condition was proved by Bryuno in \cite{Bry}, where he introduced
the so-called \emph{Bryuno condition} allowing him to still
bypass the small divisors problem in the analytic category.

Then, the problem of finding the optimal arithmetic condition for analytic linearization of analytic functions near elliptic fixed
points in dimension one was solved by Yoccoz in 
\cite{Y1,Y2}, where he proved that the Bryuno condition is in fact optimal. 
Actually, in the one dimensional case much more can be proved. Let us briefly describe it.

Given $\a\in\Bbb{R}$ irrational, let $\a_0:=\{\a\}$ be the fractional part of $\a$, 
set $\b_{-1}:=1$, and for $n\geq 1$ set $\a_{n}:=\{\su{\a_{n-1}}\}$, and $\b_n:=\a_n\dots\a_0$. 
Finally define the function
\beq{Bfunction}
\mathcal{B}(\a):=\sum_{n\geq0}\b_{n-1}\log\frac{1}{\a_n}.
\eeq
Then Yoccoz proved that, if $\{q_k\}_{k\ge 0}$ is the sequence of the denominator of the $k$-th convergent of the continuous fraction of $\a$, then $\mathcal{B}(\a)<+\infty$ if and only if
\begin{equation}\label{originale}
\sum_{k\ge0}\frac{1}{q_k}\log(q_{k+1}) < \infty,
\end{equation}
i.e.~he proved that Bryuno numbers are those for which $\mathcal{B}(\a)<+\infty$. Now, recall
%In particular, he defined the Bryuno function, that measure the size of the Siegel disk. 
%In order to better describe the result by Yoccoz, let us start by recalling 
that the Siegel disk of a complex one-dimensional dynamical system $z\mapsto f(z)$ is the maximal open 
domain $\Delta$ in which $f$ is analytically conjugate to a rotation on a disk. Recall also that, if $z_0$ is the unique fixed point of $f$ inside $\Delta$ and 
$\a\in\Bbb{R}$ is such that $e^{2\pi i\a}$ is the eigenvalue of the linearization of $f$ at $z_0$,
the conformal radius of the Siegel disk is the unique $r>0$ such that there exists a conformal diffeomorphism
$\phi: B(0,r)\rightarrow \Delta$ 
satisfying $\phi(0)=z_0$, $\phi'(0)=1$; in particular $\phi$ linearizes $f(z)$ to the rotation $R_\a(z)=(z+\a)\mod\Bbb{Z}^d$.
Then Yoccoz considered the quadratic polynomial $P_{\a}(z):=e^{2\pi i\a}z+z^2$, defined the function $r:\Bbb{R}\to\Bbb{R}_+$ as
$r(\a)=0$ if $P_{\a}(z)$ does not have a Siegel disk, and as the conformal radius of the Siegel disk of $P_\a(z)$ otherwise, and proved that
\beqano
r(\a)\geq e^{-\mathcal{B}(\a)-C_1},
\eeqano
for some $C_1>0$
thus characterizing\footnote{For rational $\a$ one has $\mathcal{B}(\a):=+\infty$.} the rotation numbers $\a$ for which $P_{\a}(z)$ admits a Siegel disk, as the Bryuno numbers.

%Then, in \cite{Y2} Yoccoz characterized the rotation numbers for which there exists a Siegel disk, which are the Bryuno numbers. 
%One of the reasons to introduce the Bryuno function comes from the fact that it gives bounds for the conformal radius.
%Indeed, in \cite{Y2} Yoccoz also showed that there exists $C>0$ such that, for all $\a\in\Bbb{R}$:\footnote{For rational $\a$, $\mathcal{B}(\a):=+\infty$.}
%\beqano
%r(\a)\geq e^{-\mathcal{B}(\a)-C}.
%\eeqano
Motivated by the relations mentioned above between the Bryuno function and the conformal radius, one may consider the function
$\Upsilon(\a):=\log r(\a)+\mathcal{B}(\a)$. Regularity properties of the function $\Upsilon$ were first guessed numerically by Marmi in \cite{M}. Later, in \cite{MMY} it was conjectured that $\Upsilon$ is $(1/2)$-H\"older continuous. The latter is the so-called Marmi-Moussa-Yoccoz conjecture:
for further developments on this topic see for instance \cite{MMY1,CC,C,Ch}.
The continuity of $\Upsilon$ (restricted the the set of Bryuno numbers) was proved by Buff and Cheritat in \cite{BC1}. 
Moreover, it was also shown by Buff and Cheritat in \cite{BC,BC1} that there exists $C_2>0$ such that for all $\a\in\Bbb{R}$ one has
\beqano
r(\a)\leq e^{-\mathcal{B}(\a)+C_2}.
\eeqano

The optimality of the Bryuno condition in the setting of analytic circle diffeomorphisms that are close to a rotation, was proved by Yoccoz in \cite{Y3}. 
Actually in \cite{Y3} much more is proved: indeed
he also provided the optimal arithmetic condition to linearize analytic circle diffeomorphisms without the assumption of being
close to rotation, completing the global theory for analytic circle diffeomorphisms that started with the work of Herman in \cite{H2}.
In the case of diffeomorphisms of the torus it is known since the work of Herman \cite{H} that such a global theory is not possible in higher dimension.
To be precise, by the so-called ``subharmonicity trick" Herman provided analytic diffeomorphisms of the two-dimensional torus with a well-defined rotation 
vector and with positive Lyapunov exponent: of course the latter implies that such a diffeomorphism cannot be analytically conjugated to a rotation. 
In particular, in this way Herman proved the ``lack of globality of KAM'' in higher dimensions. The domain of applicability of KAM theory for those diffeomorphisms 
of the torus that comes from projectivization of $SL(2,\Bbb{R})$ analytic cocycles over rotation is one of the contents of the Avila's almost 
reducibility conjecture,  that was completely solved by Avila in \cite{avilaglobal,avila_global}.
It is worth mentioning that, even if there is no global theory in higher dimension, in the case of surfaces with higher genus the situation is different and, at least 
in the case of surfaces with genus two, it is still possible to provide an analogue of Herman's global theory, as shown in work of Ghazouani-Ulcigrai in \cite{GU,GU1} (see also \cite{U} 
for a survey on diophantine-like conditions on surfaces with higher genus, that are more subtle then in the genus-one case).
In conclusion, there is ``lack of globality of KAM'' in dimension $\ge2$. 

\medskip

On the other hand, even in the local case, the problem of finding the optimal
condition on the rotation number allowing conjugacy to a roation is still open:
%As discussed above, there is "lack of globality" of KAM in higher dimensions.
%However, in the local case whether or not the Bryuno condition was optimal in higher dimensions was still open: 
the aim of the present paper is to show that the Bryuno condition is in fact not optimal in higher dimension. 
Precisely we write an arithmetic condtion (see Definition \ref{weakb} below) for which it is possible to bypass the small divisor problem in higer dimension, and we provide an
example of a vector $\a\in\Bbb{R}^2$ which is not a Bryuno vector but satisfies our condition.
%
%We will provide a counterexample in the case of analytic torus diffeomorphisms,
%but the same scheme applies also for holomorphic germs (indeed in that setting the proof is easier than that for torus diffeomorphisms).
%We will also provide an arithmetic condition that is weaker then Bryuno for which it is possible to bypass the small divisors problem in higher dimension. 
%We shall provide the KAM scheme only for the counterexample, but it works in the same way for the weaker arithmetic condition in Definition \ref{weakb}.

Let us start by introducing some notation.
Given any $v\in\Bbb{R}^d$ we denote
\begin{equation}\label{normezeta}
|v|:=\max_{j=1,\ldots,d} |v_j|,\qquad\qquad |v|_p:=\Big(\sum_{j=1}^d |\ell_j|^p\Big)^{1/p}.
\end{equation}

\begin{definition}
We denote with ${\rm{Diff}}^{\omega}_{+}(\Bbb{T}^d_{\delta})$ the group of orientation preserving homeomorphisms of $\Bbb{T}^d:=(\Bbb{R}/\Bbb{Z})^d$ which are ananalytic on the strip 
{$\Bbb{T}^d_\delta:=\{x+iy \,:\, x\in\Bbb{T}^d,\, y\in\Bbb{R}^d,  |y|<\delta\}$}.
We denote with $D^{\omega}(\Bbb{T}^d_{\delta})$ the universal covering of 
${\rm{Diff}}^{\omega}_{+}(\Bbb{T}^d_{\delta})$.

\end{definition}

\rem{nonovvio}
For $f\in D^{\omega}(\Bbb{T}^d_{\delta})$ the existence of the rotation vector i.e. of the limit 
$$
\lim_{n\to\infty} \frac{f^n}{n}({\rm mod}\Bbb{Z}^d)
$$
 is not obvious. 
However such limit is known to exist in some special cases.
Of course the limit exists if the diffeomorphism is $C^0$ conjugate to a rotation. Another case in which the limit is known to exist is the following; see for instance \cite{MJ} for details.
% that if one considers
%The existence of the rotation vector is known also in other special cases, for instance if we take 
Consider
an $SL(2,\Bbb{R})$ valued cocycle homotopic to the identity where the base is an irrational rotation on the circle, 
so that the fibered rotation number can be defined. Then, if we take the corresponding projective cocycle, it defines an orientation preserving 
homeomorphism of the torus for which there exists the rotation vector.
\erem

The classical Bryuno condition reads as follows.

\begin{definition}\label{bryvec}
For any $d\geq1$, let $\a=(\a_1,\dots,\a_d)\in\Bbb{R}^d$ a non resonant vector, i.e. such that $\a\cdot\ell\not=0$ for all $\ell\in\Bbb{Z}^d\setminus\{0\}$. For $N\in\Bbb{N}$ set
\beq{omegone}
\Omega_\a(N):=\min_{\substack{\ell\in\Bbb{Z}^d\setminus\{0\} \\ |\ell |\leq N} }\| \a\cdot\ell \|,
\eeq
where for $x\in\Bbb{R}$, we denoted
\beq{interi}
\|x\|:=\min_{l\in\Bbb{Z}}|x-l|. 
\eeq
Then $\a$ is a Bryuno vector if
\beq{bry}
\sum_{k\ge1}\frac{1}{2^k}\log\su{\Omega_\a(2^k)}<+\infty.
\eeq
\end{definition}

We shall prove the linearization Theorem \ref{main} under the following weaker condition.

\begin{definition}\label{weakb}
Let $\a=(\a_1,\dots,\a_d)\in\Bbb{R}^d$ such that $\a_1,\dots,\a_d$ are rationally independents.
 For $\ell\in\Bbb{Z}^d$, $\b\in \Bbb{S}^{d-1}:=\{y\in\Bbb{R}^d:|y|_2=1\}$ and $\delta>0$  write, 
\[
\Phi(\ell,\b,\delta):=e^{2\pi \|\ell\cdot \b\|\delta}.
\]
Fix $C>0$, $N\in\Bbb{N}$ and a sequence $c=\{c_n\}_{n\in\Bbb{N}} \in \ell^1(\Bbb{R}^+)$, and define
\beq{gic}
g(c,n,\ell):=
\left\{
\begin{aligned}
&1,\qquad \qquad2^{n}<|\ell|\leq 2^{n+1} \\
&e^{-2^n{c_n}},\qquad 0<|\ell|\leq 2^{n}, \\
&0\qquad\qquad \mbox{ otherwise}.
\end{aligned}
\right.
\eeq
Fix also $\gotd>0$ and set $\delta_{\b,n}:=\gotd$ for $n=1,\dots N$ and $\b\in\Bbb{S}^{d-1}$, while for $n\geq N$, define
\begin{equation}\label{deltini}
\delta_{\b,n+1}:=\min_{|\b-\bar{\b}|<e^{-{2^{n} c_n}}}\tilde{\delta}_{\bar{\b},n+1},
\end{equation}
where, iteratively
\begin{equation}\label{tildenpiuno}
\begin{aligned}
\tilde{\delta}_{\b,n+1}:=\max\{&\delta\leq \delta_{\b,n} \ : \ \frac{\Phi(\ell,\b,\delta)}{\|\a\cdot\ell\|} g(c,n,\ell)  \leq \max_{\bar{\b}\in\Bbb{S}^{d-1}}\Phi(\ell,\bar{\b},\delta_{\bar{\b},n}) \\
&\quad\forall\ 0< |\ell|\leq 2^{n+1}\}.
\end{aligned}
\end{equation}
We say that $\a$ is a \emph{weak-Bryuno vector} if for every $\gotd>0$ there exists a choice of the constants $C$, 
the sequence $c=\{c_n\}_{n\in\Bbb{N}}$ and the step $N\in\Bbb{N}$ %and $\delta_{\b,n}$ 
as above such that
\beq{weakcond}
\inf_{\b\in\Bbb{S}^{d-1}}\inf_{n\in\Bbb{N}}\delta_{\b,n}>0.
\eeq
\end{definition}

Roughly, for some choice of $c$, the weight $g(c,n,\ell)$ is small enough to control the small divisors for $|\ell|<2^n$,
while the sequence $\{\delta_{\b,n}\}_{n\ge1}$ is chosen in order to to solve the cohomological equation for $2^{n}<|\ell|\leq 2^{n+1}$.
In Appendix \ref{appendice} we prove that a number is weak-Bryuno if and only if it is a Bryuno number, thus implying that the difference occurs only in
dimension $\geq2$.

Our first result, proved in Section \ref{superprova}, is the following.

\thm{super}
There exists $\a\in\Bbb{R}^2$  that is weak-Bryuno vector but not a Bryuno vector.
\ethm

\rem{commento1}
Let us comment the arithmetic condition above. 
In order to prove a linearization result with the condition  in Definition \ref{weakb}, at each KAM step
% we cannot use the standard KAM scheme,  but 
we need to loose regularity in a non uniform way. Precisely, at each scale $2^n$ we divide the fundamental domain of the torus into
slices, and in each slice we will have an appropriate loss of regularity
%Instead of the usual KAM schemes, in order to prove linearization results with the condition 
%\eqref{weakcond} we will loose regularity not in a uniform way but on each scale we will divide the fundamental domain of the Torus in slices and 
%in each of those slices we will define an appropriate lost of regularity 
(which will be  $\sim(\delta_{\b,n}-\delta_{\b,n+1})$, with $\b$ parametrizing the direction).
The r\^ole of the sequence $\{c_n\}_{n\in\Bbb{N}}$ is the following: at each KAM step we conjugate the diffeomorphism 
with a map close to the identity in such a way that the new reminder is essentially quadratically smaller. Then, the r\^ole of $\{c_n\}_{n\in\Bbb{N}}$ 
is to ensure that the slices at scale $2^n$ are mapped via the conjugacy into slices at scale $2^{n+1}$; in other words we need to make sure that at each step 
we do not simply have loss of regularity of the thickened torus, but also the slices have to be enlarged a bit.
\erem

\rem{commento2}

It is not difficult to show that the scheme converges also if one sets $g(c,n,\ell)=0$ for $0<|\ell|\leq 2^n$. Indeed after solving the cohomological 
equation at scale $2^n$, one can always repeat the argument and solve the cohomological equation with truncation $2^n$ a finite number of times
 in order  to make the truncation at scale $2^n$ of the error term as small as one whishes.
\erem

For any $f\in {\rm D}^{\omega}_{+}(\Bbb{T}^d_{\delta})$ introduce the norm
\beq{normana}
\|f \|_{\delta}:=\sum_{\ell\in\Bbb{Z}^2} |\hat{f}(\ell)| e^{2\pi \delta|\ell|_2},
\eeq
and for any $\a\in\Bbb{R}^d$ let $R_\a$ denote the rotation with rotation vector $\a$, i.e. $R_\a(z) =(z+\a){\rm mod}\Bbb{Z}^d$.
We shall prove the following result.

\thm{main}
Fix any weak-Bryuno vector $\a\in\Bbb{R}^d$ and  $\delta>0$. Then there is $\e=\e(\a,\delta)>0$ such that, for any $f\in {\rm D}^{\omega}(\Bbb{T}^d_{\delta})$
 such that
$$
\lim_{n\to\infty}\frac{f^n}{n}({\rm{mod}}\Bbb{Z}^d)=\a ({\rm{mod}}\Bbb{Z}^d)
$$ 
and $\|f-R_{\a}\|_{\delta}<\e$, there exists 
$H\in \Diff^{\omega}(\Bbb{T}^d_{\frac{\delta}{2}})$ with $\|H-id.\|_{\frac{\delta}{2}}<\sqrt{\e}$ such that
\beq{}
H^{-1}\circ f\circ H=R_\a.
\eeq
\ethm

Let us add some remarks.

\medskip

%\giu
%\nl
%{\bf{Optimal arithmetic condition in higher dimension}}  We expect the condition \eqref{weakb} to be optimal. Indeed blabla
%\red{sicuri che vogliamo metterlo???}

\giu
\nl
{\bf{Similar results for other type of systems}}
We believe that
a result similar to Theorem \ref{main} can be obtained for other types of dynamical systems. For instance we expect that, by
using a scheme similar to the one in the present paper,
it would be possible to show that the weak-Bryuno condition of Definition \ref{weakb} is enough to obtain
the stability of analytic area-preserving mapping near a fixed point in $\Bbb{R}^{d}$ for $d\geq 3$. In particular, this would be an improvement of \cite{Rus1}. 
In the same way, a similar statement could also be obtained for instance for holomorphic germs.

\giu
\nl
{\bf{Optimal arithmetic condition for non-analytic classes }}
A corresponding Bryuno-arithmetic condition was defined also in other contexts. For instance, for $\theta\in\Bbb{R}^d$,
in ref.~\cite{BF}, the so-called $\alpha$-Bryuno-R\"ussmann condition, i.e.~the condition that
\begin{equation}\label{buofe}
\mathcal{F}_{\alpha}(\theta):=\int_1^{+\infty} \frac{1}{Q^{1+\frac{1}{\alpha}}}\log\left(\frac{1}{\beta_\theta(Q)}\right) d Q <\infty,
\qquad \alpha\ge 1,
\end{equation}
is used to obtain KAM tori in the Gevrey class.  Note that the set $BR_\alpha:=\{\theta\in\Bbb{R}^d:\mathcal{F}_\alpha(\theta)<\infty\}$ decreases
w.r.t.~$\alpha$, and for $\alpha=1$ one obtains KAM tori in the analytic case, as shown by R\"ussmann in ref.~\cite{Rus}.
Note that the condition \eqref{originale} is equivalent to \eqref{buofe} with $\alpha=1$; see also
ref.~\cite{Bounemoura}, where these and other related issues are discussed. In \cite{Bounemoura 1} Bounemoura shows that this condition is not optimal and 
he asks whether for some quasi-analytic classes the corresponding Bryuno condition is optimal. 
We think that a method similar to the one in the present paper could be implemented to show that
also in case of quasi-analytic classes the corresponding Bryuno-type condition is not optimal.

\giu
\nl
{\bf{Infinite dimensional case }}
Recently an infinite dimensional version of the Bryuno condition was introduced in \cite{CGP} to obtain the exsitence of almost-periodic
solutions for a family of NLS equations with smooth convolution potential. It would be interesting to study whether the ideas of the present
paper can be applied also to the infinite-dimensional case.

\giu
\nl
{\bf{Further results}} In \cite{KQY}, Krikorian, Wang, You and Zhou proved the non optimality of the Bryuno condition in higher dimension 
in the non-linear skew product setting. However our result is quite different and holds for any small analytic perturbation of a rotation for which the rotation vector exists.
\medskip

The paper is organized as follows. In Section \ref{superprova} we construct a rotation vector $\a=(\a_1,\a_2)\in\Bbb{R}^2$ that is weak-Bryuno but not Bryuno. 
In particular we write the components $\a_1,\a_2$ as series whose summands decrease very fast\footnote{ 
This is somehow reminiscent of the method used by Liouville to construct the first examples of transcendental numbers.} and an asymptotic behavior which
is essentially the same. 
In other words, if we write 
$\a_1=\sum_{n\ge0}\su{b_{n,1}}$, $\a_2=\sum_{n\ge0}\su{b_{n,2}}$, then $b_{n,1}\sim b_{n,2}$ while   $b_{n+1,i}\ll b_{n,i}$. Moreover
we need that each scale $2^n$ there is $\ell_n\in\Bbb{Z}^2$ such that $\|\ell_n\cdot\a\|\sim b_{n+1,1}\sim b_{n+1,2}$,
%of $\a_1,\a_2$ such that that the maximal loss of regularity is comparable to the next term of the series of $\a_1,\a_2$, and so that at different 
so that at different scales scales the direction with maximal loss of regularity is twisted. 
We write the components $\a_1,\a_2$ as series in order to clearly identify the derection of maximal loss of regularity at each scale and twist it at the next scale.
%The fact that $\a_1,\a_2$ are not written in terms of their continued fraction  but as series helps to find these linear combinations and to twist the direction. 

In Section \ref{provakam} we prove Theorem \ref{main}.
%i.e.~we provide a KAM scheme to obtain the 
%map that analytically linearizes analytic diffeomorphisms of the torus with rotation vector $\a$ and that are close to $R_{\a}$. 
In particular, in order to 
prove the convergence of the scheme, we introduce some new norms (see \eqref{stanormina}) that detect the fact that at each scale there 
is essentially only one direction where there is some significant loss of regularity. In particular, at each scale instead of a uniform loss of 
regularity, we divide a fundamental domain of the torus into slices of size related to the maximal loss of regularity, and in different 
slices we obtain a different loss of regularity. Finally, before going to the subsequent scale, we repeat the scheme a few times by solving the 
cohomological equation with the same truncation, in order to get the truncation of the error term arbitrarily small; see also Remark \ref{commento2}.

\medskip

\noindent
{{\bf Acknowledgements}}. We  thank A.~Avila, B.~Fayad, G.~Gentile and M.~Procesi for many useful discussion
and suggestions. We also thank an anonymous referee for encouraging us to improve the first draft of the manuscript.
F.A.~has been supported by the SNSF grant.
L.C.~has been  supported by the  research projects  PRIN 2020XBFL ``Hamiltonian and dispersive PDEs" 
and PRIN 2022HSSYPN ``Turbulent Effects vs Stability in Equations from Oceanography'' (TESEO)
of the 
Italian Ministry of Education and Research (MIUR).

%%%%%%%%%%%%%%%%%%%%%%%%%%%%%%%%%%%%%%%%%%%%%%%%%%%%%%%%%%%%%%%%%%%%%%%%%%%%%%%%

\zerarcounters
\section{A class of weak-Bryuno vector which are not Bryuno vectors}\label{superprova}
As already mentioned,
in dimension 1 the Bryuno condition is equivalent to (\ref{originale}); in particular, in order to obtain a linearization result one has to control the loss of regularity at scale $q_n$.
In dimension 2 we construct a vector $\a=(\a_1,\a_2)$ so that
both $\a_1,\a_2$ satisfy the Bryuno condition (\ref{originale}), 
and, similarly to the one-dimensional case, there is a sequence $\{\tilde{q}_l\}_{l\in\Bbb{N}}$
 such that 
 
 \begin{enumerate}

\item 
 at each scale $\tilde{q}_l$ there exists only one direction $\nu(l)\in\Bbb{Z}^2$ for which $\|\a\cdot\nu(l)\|$ is really small,
 
 \item
 the ``smallest direction'' is twisted enough when we consider the subsequent scale $\tilde{q}_{l+1}$.
 
 \end{enumerate}

% , whereas in the other directions there is much 
% less loss of regularity. Moreover, as already mentioned, since in the proof of Theorem \ref{main}
% 
% Moreover, in order to prove convergence of the scheme we will need to look at the loss of regularity in each direction on these scales.  
% Finally, we want also that this ``worst" direction to be twisted when we change the scale. 

Let us first introduce the following auxiliary sequence.
 
\begin{definition}\label{ricordatela}
Fix $0<\bar{\theta}<\su{2}$. We define the sequence $\{t_n\}_{n\in\Bbb{N}}$ recursively as
\beq{tn}
t_0=2,\qquad  t_{n+1}:=t_n+\left\lfloor t_n^{\bar{\theta}}\right
\rfloor.
\eeq
\end{definition}

\rem{tienne}
By construction one has $t_n\ge n$. More precisely
\[
t_{2n}\ge n + \sum_{j=n}^{2n}j^{\bar{\theta}}\ge n^{1+\bar{\theta}},
\]
and in particular
\[
t_n \ge Cn^{1+\bar{\theta}},
\]
for some constant $C=C(\bar{\theta})$.
\erem

As we shall see, the r\^ole of the sequence $\{t_n\}_{n\ge0}$ is the following: if at the $t_n$-th step of the KAM iteration the
``worst'' loss of regularity is in the direction $v$, then $t_{n+1}$ is the first step at which the ``worst'' loss of regularity is in a
direction $1/t_{n+1}^{\frac{\bar{\theta}}{2}}$-close to that of $v$.

The aim of this section is to prove the Proposition below which provides a class of vectors which are weak-Bryuno
vectors but not Bryuno vectors. Then, in the subsequent section we will show that such a class is non-empty.

\pro{criterio}
Let $(\a_1,\a_2)\in\Bbb{R}^2$ and fix $\bar{\theta}\in(0,\su{2})$ and $\{t_n(\bar{\theta})\}_{n\in\Bbb{N}}$ as in Definition \ref{tn}. 
Suppose that there exist $\{\tilde{q}_l\}_{l\in\Bbb{N}}\subset\Bbb{N}$ increasing and $\{\tilde{\nu}(l)\}_{l\in\Bbb{N}}\subset\Bbb{Z}^2$ such that
the following holds.

\begin{enumerate}

\item $\tilde{q}_{l+1}>4\tilde{q}_l$ for every $l\in\Bbb{N}$.

    \item For every $l\in\Bbb{N}$ one has $|\tilde{\nu}(l)|=\tilde{q}_l$, $\Omega_{\a}(\tilde{q}_l) = \frac{1}{\|\a \cdot \tilde{\nu}(l)\|}$ and
    \[
   \frac{1}{2t_n\log t_n}\leq \frac{1}{\tilde{q}_l}\log\su{\Omega_{\a}(\tilde{q}_l)}\leq\frac{4}{t_n\log t_n} ,
    \]
    with $n=n(l)$ such that $t_n<l\leq t_{n+1}$, for some $C>0$.

    \item For every $l_1,l_2\in\Bbb{N}$ such that $t_n<l_1<l_2\leq t_{n+1}$ there is $C>0$ such that
    \[
    \theta(\tilde{\nu}(l_1),\tilde{\nu}(l_2))\geq {C}{t_n^{-\frac{\bar{\theta}}{2}}},
    \]
    where $ \theta(\tilde{\nu}(l_1),\tilde{\nu}(l_2))$ denotes the angle between $\tilde{\nu}(l_1)$ and $\tilde{\nu}(l_2)$.

    \item For $n\in\Bbb{N}$, setting $l(n):=\max \{l\in\Bbb{N}:\tilde{q}_l\leq 2^n\}$ and
    \[
    k_n:=\min_{\substack{ \ell\in\Bbb{Z}^2 \\ \ell \nparallel \nu(l(n)) \\ |\ell|\le 2^n}} \| \a \cdot\ell \|,
    \]
    Then
    \[
    \sum_{n\in\Bbb{N}}\frac{1}{2^n}\log\su{k_n} <\infty.
    \]
    
    \item
    Let $\ell\parallel \tilde{\nu}(l)$.
    For $\ell\in\Bbb{Z}^2$ such that $\ell\parallel \tilde{\nu}(l)$, if $\tilde{q}_l\leq|\ell|\leq \sqrt{\tilde{q}_{l+1}}$ we have
    \beq{ncora}
    \frac{1}{|\ell|}\log \su{\|\a\cdot\ell\|} \leq \frac{C\tilde{q}_l}{2^ n t_m\log t_m},
    \eeq
    while if $\sqrt{\tilde{q}_{l+1}}<|\ell|<\tilde{q}_{l+1}$ we have
    \beq{ncora1}
    \frac{1}{2^n}\log \su{k_n} \leq \frac{C\sqrt{\tilde{q}_{l+1}}}{2^n t_m\log t_m},
    \eeq
    
\end{enumerate}

Then, $\a$ is a weak-Bryuno vector but it is not a Bryuno vector.
\epro

Before proving Proposition \ref{criterio} let us comment the conditions 1--5.
The first two conditions guarantee that the vector is not a Bryuno vector. Indeed, as we shall see, %\footnote{Below we will give more details about the following inequality.}
it implies
\beqano
\sum_{l\in\Bbb{N}}\frac{1}{\tilde{q}_l} \log\su{\Omega_{\a}(\tilde{q}_l)}\geq C\sum_{n\in\Bbb{N}}\frac{1}{\log t_n}=\infty.
\eeqano

The fourth condition means that all vectors at scale $\tilde{q}_l$ which are not parallel to $\tilde{\nu}(l)$ are negligible, in the sense that the corresponding small divisors
are not that small compared to $\langle\a,\tilde{\nu}(l)\rangle$. 
The other conditions guarantee that $\a$ is a weak-Bryuno vector: indeed as we shall see, if the angle between a fixed direction $v$ and the ``worst direction''
 $\frac{\tilde{\nu}(l)}{|\tilde{\nu}(l)|_2}$ is big enough, the loss of regularity along the direction $v$ will be of the 
same order as the loss of regularity in the direction of integer vectors not parallel to $\tilde{\nu}(l)$,
%given by integer vectors not parallel to the worst direction, 
which is controlled by condition 2.

Let us first show that if $\a$ satisfies the conditions of Proposition \ref{criterio}, then it is not a Bryuno vector.

\lem{}
Let $\a\in\Bbb{R}^2$ satisfying conditions 1--5 of Proposition \ref{criterio}. Then, $\a$ is not a Bryuno vector.
\elem

\proof
For $l\in \Bbb{N}$ let $m(l)\in\Bbb{N}$ be such that $2^{m(l)-1}<\tilde{q}_l\leq 2^{m(l)}$. 
Then one has
\[
\begin{aligned}
\sum_{n\in\Bbb{N}}\frac{1}{2^n}\log\su{\Omega_{\a}(2^n)}&\geq \sum_{l\in\Bbb{N}}\frac{1}{2^{m(l)}}\log\su{\Omega_{\a}(2^{m(l))})}\\
&\geq \sum_{l\in\Bbb{N}}\frac{1}{2\tilde{q}_l}\log\su{\Omega_{\a}(\tilde{q}_l)}\geq C\sum_{l\in\Bbb{N}}\su{l\log l}=+\infty,
\end{aligned}
\]
where in the last inequality we have used the fact that for $t_n<l\leq t_{n+1}$ it holds
\[
\frac{1}{2\tilde{q}_l}\log\su{\Omega_{\a}(\tilde{q}_l)}\geq \frac{C}{t_n\log t_n}\geq \frac{C}{l\log l}.
\]
In particular, $\a$ is not a Bryuno vector. 
\qed

We are left to prove that conditions 1--5 imply that $\a$ is a weak-Bryuno vector.
Let us start by fixing some notations.

\begin{definition}\label{cappa}
    For $l\in \Bbb{N}$ we define $m(l)\in\Bbb{N}$ as the unique integer such that
    \[
    2^{m(l)-1}<\tilde{q}_l\leq 2^{m(l)}.
    \]
    For $n\in\Bbb{N}$ we define $l(n)$ as the unique integer such that
    \[
    \tilde{q}_{l(n)}\leq 2^{n}<\tilde{q}_{l(n)+1}
    \]
    and $k(l)$ as the unique integer such that
    \[
    t_{k(l)}<l\leq t_{k(l)+1} .
    \]
\end{definition}

\pro{altra}
Suppose that $\a\in\Bbb{R}^2$ satisfies the conditions of Proposition \ref{criterio}. Then, $\a$ is a weak-Bryuno vector. 
\epro

We need to show that for every $\gotd>0$ there exists a choice of $C>0, N\in\Bbb{N}$ and $\{c_n\}_{n\in\Bbb{N}}$ such that the sequences
$\{\delta_{\b,n}\}_{n\in\Bbb{N}}$ satisfy
\[
\inf_{\b\in\Bbb{S}^{d-1}}\inf_{n\in\Bbb{N}}\delta_{\b,n}>0.
\]

Set
\[
c_n:=\frac{C}{t_{k(l(n))}(4n-\log\tilde{q}_{l(n)})^2}+\frac{C}{t_{k(l(n))}(|4n-\log\tilde{q}_{l(n)+1}|+1)^2}+\frac{1}{2^n}\log \su{k_n},
\]
with $C,t_n, k_n$ so that the conditions in Proposition \ref{criterio} hold.

\lem{log}
There exists $\t>1$ such that, for all $n\in\Bbb{N},\ell\in\Bbb{Z}^2$ with $|\ell|\leq 2^n$ and $2^{n+(\log t_{k(l(n))})^{\t})}<\tilde{q}_{l(n)+1}$ we have
\[
g(c,\ell,n+(\log t_{k(l(n))})^{\t})\leq e^{-2^n \frac{t_{k(l(n))}^2}{(4n-\log\tilde{q}_{l(n)})^2}}e^{-2^n \frac{t_{k(l(n))}^2}{t_{k(l(n))}(|4n-\log\tilde{q}_{l(n)+1}|+1)^2}}.
\]
\elem
\proof
By the definition of $g$ and the fact that $|\ell|\leq 2^n$ we obtain
\[
\begin{aligned}
g(c,\ell,n+(\log t_{k(l(n))})^{\t})&=e^{-2^{n+(\log t_{k(l(n))})^{\t}}c_{n+(\log t_{k(l(n))})^{\t}}} \\
&\leq e^{-\frac{2^{n+(\log t_{k(l(n))})^{\t}}}{t_{k(l(n))}(4n-\log\tilde{q}_{l(n)})^2(\log t_{k(l(n))})^{\t}}}e^{-\frac{2^{n+(\log t_{k(l(n))})^{\t}}}{t_{k(l(n))}(4n-\log\tilde{q}_{l(n)})^2(\log t_{k(l(n))})^{\t}}} \\
& \leq e^{-2^n \frac{t_{k(l(n))}^2}{t_{k(l(n))}(|4n-\log\tilde{q}_{l(n)+1}|+1)^2}}e^{-2^n \frac{t_{k(l(n))}^2}{t_{k(l(n))}(|4n-\log\tilde{q}_{l(n)+1}|+1)^2}}
\end{aligned}
\]
for $\t$ large enough.
\qed

\giu
\nl

%As a corollary of Lemma \ref{log} together with the bounds \eqref{ncora}, \eqref{ncora1} we obtain the following result.
\cor{}
Let $\ell\in\Bbb{Z}^2$, $n\in\Bbb{N}$ be such that $|\ell|\leq 2^{n}$ and $\ell\parallel \tilde{\nu}(l(n))$. Then, for $p\geq (\log t_{k(l(n))})^{\t} $ such that 
$2^{n+p}< \tilde{q}_{l(n)+1}$ we have
\[
\frac{g(c,\ell,n+p)}{\|\a\cdot\ell\|}\leq 1.
\]
\ecor
\proof
It follows directly from Lemma \ref{log} using the bounds \eqref{ncora}, \eqref{ncora1}.
\qed

\begin{definition}
For $l\in\Bbb{N}$ let $n(l)\in\Bbb{N}$ such that $2^{n(l)}<\tilde{q}_{l+1}\leq 2^{n(l)+1}$  and let $\b_l:=\frac{\tilde{\nu}(l)}{|\tilde{\nu}(l)|_2}$, with $\tilde{\nu}(l)$ as in Proposition \ref{criterio}. 
\end{definition}

The fact that $\a$ is a weak-Bryuno vector will follow from the Proposition below.
\pro{ttn}
Let $k\in\Bbb{N}$, $\delta>0$ and suppose that:
\[
\inf_{\b\in\Bbb{S}^{d-1}}\delta_{\b,n(t_k)}\geq \delta.
\]
Then one has
\[
\inf_{\b\in\Bbb{S}^{d-1}}\delta_{\b,n(t_{k+1})}\geq \delta-\frac{C(\log t_{k})^{\t}}{t_k}-C \sum_{\tilde{q}_{t_{k+1}}\leq 2^j<\tilde{q}_{t_{k+2}}}c_j.
\]
\epro
Let us first show that the latter result implies that $\a$ is a weak-Bryuno vector. 

\giu
\nl
{\bf{Proof (Proposition \ref{ttn} implies Proposition \ref{altra})}}
By definition of the sequence $\{\delta_{\b,n}\}_{n\in\Bbb{N},\b\in\Bbb{S}^{d-1}}$ we know that for all $\b\in\Bbb{S}^{d-1},n\in\Bbb{N}$ we have
\[
\delta_{\b,n}\geq \delta_{\b,n+1}.
\]
Indeed by \eqref{tildenpiuno} one has
\[
\tilde{\delta}_{\b,n+1}\leq \delta_{\b,n},
\]
which in turn by \eqref{deltini} implies
\begin{equation}\label{delti}
\delta_{\b,n+1}\leq \tilde{\delta}_{\b,n+1}\leq\delta_{\b,n}.
\end{equation}
In particular, %, if $N$ is such that $\delta_{\b,n}=1$ for $1\leq n\leq N,\b\in\Bbb{S}^{d-1}$, we have
from the fact that $t_n\geq Cn^{1+\bar{\theta}}$ and the fact that $\{c_n\}_{n\in\Bbb{N}}\in \ell^1(\Bbb{R}^{+})$, we deduce
\[
\inf_{\b\in\Bbb{S}^1,n\in\Bbb{N}}\delta_{\b,n}=\inf_{\b\in\Bbb{S}^1,n\in\Bbb{N}}\delta_{\b,t_n}\geq 1-C\sum_{n\geq N}\left(\frac{(\log t_{n})^c}{t_n}+c_n\right)>0
\]
if $N$ is large enough. % Indeed, from the fact that $t_n\geq Cn^{1+\bar{\theta}}$ and the fact that $\{c_n\}_{n\in\Bbb{N}}\in \ell^1(\Bbb{R}^{+})$ we know that the sum is convergent.
\qed

\giu
\nl

We now show that, chosing $\{c_n\}_{n\ge1}$ above, the function $g(c,\ell,n)$ in Definition \ref{weakb} is much smaller then the small divisor
$\|\a\cdot \ell\|$ for $|\ell|\leq 2^{n}$ and $\ell\nparallel \nu(l(n))$.
\lem{x}
Let $n,\ell\in\Bbb{N}$ such that $|\ell| \leq 2^{n}$ and $\ell\nparallel \nu(l(n))$. Then
\[
\frac{g(c,\ell,n)}{\|\a\cdot\ell\|}\leq 1.
\]
\elem
\proof
Recall that for $n\in\Bbb{N}$ by definition of $c_n$ we have that $c_n\geq\frac{1}{2^n}\log\su{k_n}$. Then, by definition of $g(c,n,\ell)$ and for $|\ell|\leq 2^n$
we have
\[
g(c,n,\ell)=e^{-2^n c_n}\leq %e^{-2^n(\frac{\log\su{k_n}}{2^n})}=
\su{k_n}\leq \|\a\cdot\ell\|,
\]
with the last inequality following from the definition of $k_n$ and the fourth condition in Proposition \ref{criterio}.
\qed

\giu
\nl

The result above essentially implies that the only small divisors we need to worry about are the ones which are parallel to $\tilde{\nu}(l)$.

\giu
\nl

Let us now provide a basic bound for the sequence $\{\delta_{\b,n}\}_{n\ge1}$.
\lem{facciamolo}
Fix $n\in\Bbb{N}$, $\{\delta_{\b,n}\}_{\b\in\Bbb{S}^1}$, $\{\bar{\delta}_{\b,n}\}_{\b\in\Bbb{S}^1}$  and suppose that for all $\b\in \Bbb{S}^1$ one has
\[
\delta_{\b,n}\geq \bar{\delta}_{\b,n}.
\]
Define
\begin{equation}\label{deltahat}
\begin{aligned}
\hat{{\delta}}_{\b,n+1}:=\max\{&\delta\leq \bar{\delta}_{\b,n} \ : \ \frac{\Phi(\ell,\b,\delta)}{\|\a\cdot\ell\|} g(c,n,\ell)  \leq \max_{\bar{\b}\in\Bbb{S}^{d-1}}\Phi(\ell,\bar{\b},\bar{\delta}_{\bar{\b},n}) \\
&\quad\forall\ 0< |\ell|\leq 2^{n+1}\},
\end{aligned}
\end{equation}
and set

\begin{equation}\label{delitti}
\bar{\delta}_{\b,n+1}:=\min_{|\b-\bar{\b}|<e^{-{2^{n} c_n}}}\hat{{\delta}}_{\bar{\b},n+1}.
\end{equation}

Define also
\begin{equation}\label{deltabar}
\begin{aligned}
\tilde{\delta}_{\b,n+1}:=\max\{&\delta\leq \delta_{\b,n} \ : \ \frac{\Phi(\ell,\b,\delta)}{\|\a\cdot\ell\|} g(c,n,\ell)  \leq \max_{\bar{\b}\in\Bbb{S}^{d-1}}\Phi(\ell,\bar{\b},\delta_{\bar{\b},n}) \\
&\quad\forall\ 0< |\ell|\leq 2^{n+1}\},
\end{aligned}
\end{equation}
and set
\begin{equation}\label{deltizzi}
\delta_{\b,n+1}:=\min_{|\b-\bar{\b}|<e^{-{2^{n} c_n}}}\tilde{\delta}_{\bar{\b},n+1}.
\end{equation}

Then, for all $\b\in\Bbb{S}^1$ we have
\[
\delta_{\b,n+1}\geq\bar{\delta}_{\b,n+1}.
\]

\elem
\proof
Since for all $\b\in \Bbb{S}^1$ we have $\bar{\delta}_{\b,n}\leq \delta_{\b,n}$, then, for all $\ell\in\Bbb{Z}^2$ such that $|\ell|<2^{n+1}$ we have:
\[
\max_{\b\in\Bbb{S}^1}\Phi(\ell,\b,\bar{\delta}_{\b,n})\leq \max_{\b\in\Bbb{S}^1}\Phi(\ell,\b,\delta_{\b,n}),
\]
from which we deduce
\[
\hat{{\delta}}_{\b,n+1}\leq \tilde{\delta}_{\b,n+1}.
\]
But then
\[
\bar{\delta}_{\b,n+1}:=\min_{|\bar{\b}-\b|\leq e^{-2^n c_n}}\hat{{\delta}}_{\b,n+1}\leq \min_{|\bar{\b}-\b|\leq e^{-2^n c_n}}\bar{\delta}_{\b,n+1}=\delta_{\b,n+1}
\]
so the assertion holds.
\qed

\giu
\nl

Fix $C>0$ large enough, let $C_0=C$, for $n\in\Bbb{N}$ set 
\begin{equation}\label{cione}
C_{n+1}:=C_n(1+e^{-2^n c_n}).
\end{equation}
and note that $\{C_n\}_{n\ge1}\in\ell^{\infty}(\Bbb{R}^+)$.
% In particular, note that there exists $\bar{C}>0$ such that $C_n<\bar{C}$ for all $n\in\Bbb{N}$. Moreover, 
Also, given any $l\in\Bbb{N}$, for $n\in\Bbb{N}$ such that $\tilde{q}_l\leq 2^n<\tilde{q}_{l+1}$ set
\begin{equation}\label{pienne}
p_n := 
\left\{
\begin{aligned}
&1,\qquad  &2^n<\sqrt{\tilde{q}_{l+1}}, \\
&2^{-n(t_k)}\sqrt{\tilde{q}_{l+1}},\qquad  &2^n\geq \sqrt{\tilde{q}_{l+1}},
\end{aligned}
\right.
\end{equation}
where $k=k(l)$ is as in Definition \ref{cappa}.
% $p_n=1$ if $2^n<\sqrt{\tilde{q}_{l+1}}$ and $p_n=2^{-n(t_k)}\sqrt{\tilde{q}_{l+1}}$ if $2^n\geq \sqrt{\tilde{q}_{l+1}}$ and $k=k(l)$.
Then Proposition \ref{ttn} follows directly from the following result.

\lem{bounds}
Let $k,l\in \Bbb{N}$ be such that $t_k\leq l<t_{k+1}$, with $t_k$ as in Definition \ref{ricordatela}.
Suppose that for all $\b\in\Bbb{S}^1$ the following holds.

\begin{itemize}

\item If there exists $h=t_k,\ldots, l$ for which
$$
\theta(\b,\tilde{\nu}(h))<{C_{n(l)}}{t_k^{-\frac{\bar{\theta}}{2}}}
$$
then
\[
\delta_{\b,n(l)}\geq \delta-C\sum_{n=n(t_k)}^{ n(l)}\Big(\frac{(\log t_k)^{\t}}{2^{n-n(t_k)}t_k}-c_n \Big).
\]

\item If 
$$
\theta(\b,\tilde{\nu}(h)) > {C_{n(l)}}{t_k^{-\frac{\bar{\theta}}{2}}}
$$
for all $h=t_k,\ldots, l$, then
\[
\delta_{\b,n(l)}\geq \delta-C\sum_{n=n(t_k)}^{n(l)}c_n.
\]

\end{itemize}

%
%
%
%
% such that there exists $t_k\leq h\leq l$ for which $\theta(\b,\tilde{\nu(h)})<{C_{n(l)}}{t_k^{-\frac{\bar{\theta}}{2}}}$, we have:
%\[
%\delta_{n(l),\b}\geq \delta-\sum_{n(t_k)\leq n\leq n(l)}\frac{C(\log t_k)^{\t}}{2^{n-n(t_k)}t_k}-C\sum_{n(t_k)\leq n\leq n(l)}c_n
%\]
%and for $\b\in\Bbb{S}^1$ such that it does not exists $t_k\leq h\leq l$ such that $\theta(\b,\tilde{\nu(h)})<\frac{C_{n(l)}}{t_k^{\frac{\bar{\theta}}{2}}}$, we have:
%\[
%\delta_{n(l),\b}\geq \delta-C\sum_{n(t_k)\leq n\leq n(l)}c_n.
%\]
Then, for all $\b\in\Bbb{S}^1$  the following holds.

\begin{enumerate}

\item If there exists $h=t_k,\ldots, l+1$ for which
$$
\theta(\b,\tilde{\nu}(h))<{C_{n(l+1)}}{t_k^{-\frac{\bar{\theta}}{2}}}
$$
then
\[
\delta_{\b,n(l+1)}\geq \delta-C\sum_{n=n(t_k)}^{ n(l+1)}\Big(\frac{(\log t_k)^{\t}}{2^{n-n(t_k)}t_k}-c_n \Big)
\]

\item If 
$$
\theta(\b,\tilde{\nu}(h)) > {C_{n(l+1)}}{t_k^{-\frac{\bar{\theta}}{2}}}
$$
for all $h=t_k,\ldots, l+1$, then
\[
\delta_{\b,n(l+1)}\geq \delta-C\sum_{n=n(t_k)}^{n(l+1)}c_n.
\]

\end{enumerate}
%
%
%
%
%
%such that there exists $t_k\leq h\leq l+1$ such that $\theta(\b,\tilde{\nu(h)})<\frac{C_{n(l+1)}}{t_k^{\frac{\bar{\theta}}{2}}}$, we have:
%\[
%\delta_{n(l+1),\b}\geq \delta-\sum_{n(t_k)\leq n\leq n(l+1)}\frac{C(\log t_k)^{\t}}{2^{n-n(t_k)}t_k}-C\sum_{n(t_k)\leq n\leq n(l+1)}c_n
%\]
%and for $\b\in\Bbb{S}^1$ such that it does not exists $t_k\leq h\leq l+1$ such that $\theta(\b,\tilde{\nu(h)})<\frac{C_{n(l+1)}}{t_k^{\frac{\bar{\theta}}{2}}}$, we have:
%\[
%\delta_{n(l+1),\b}\geq \delta- C\sum_{n(t_k)\leq n\leq n(l+1)}c_n.
%\]
\elem

\proof
For all $n=n(l),\ldots,n(l+1)$ and all $\b \in \Bbb{S}^1$ define 
\begin{equation}\nonumber
\tilde{\delta}_{\b,n} :=
\left\{
\begin{aligned}
&\delta-C\sum_{m=n(t_k)}^{ n}\Big(\frac{(\log t_k)^{\t}}{2^{m-n(t_k)}t_k}-c_m \Big), & \b \mbox{ is as in case 1} \\
& \delta-C\sum_{m=n(t_k)}^{n}c_m, & \b \mbox{ is as in case 2}
\end{aligned}
\right.
\end{equation}

%
%
%
%Let us define for $n(l)\leq n\leq n(l+1)$, $\b\in\Bbb{S}^1$ such that there exists $t_k\leq h\leq l+1$ such that $\theta(\b,\tilde{\nu(h)})<\frac{C}{t_k^{\frac{\bar{\theta}}{2}}}$:
%\[
%\tilde{\delta}_{n,\b}= \delta-\sum_{n(t_k)\leq m\leq n}\frac{C(\log t_k)^{\t}p_m}{2^{m-t_k}t_k}-C\sum_{n(t_k)\leq m\leq n}c_m
%\]
%and for $\b\in\Bbb{S}^1$ such that it does not exists $t_n\leq h\leq l+1$ such that $\theta(\b,\tilde{\nu(h)})<\frac{C}{t_k^{\frac{\bar{\theta}}{2}}}$:
%\[
%\tilde{\delta}_{n,\b}= \delta-C\sum_{n(t_k)\leq m\leq n}c_m.
%\]

We shall prove that $\delta_{\b,n}\geq \tilde{\delta}_{\b,n}$ implies that $\delta_{\b,n+1}\geq \tilde{\delta}_{\b,n+1}$.
Let us define also
\begin{equation}\label{deltiagain}
\bar{\delta}_{\b,n+1}:=\min_{|\b-\bar{\b}|<e^{-{2^{n} c_n}}}{\check{\delta}}_{\bar{\b},n+1},
\end{equation}
where, iteratively
\[
\begin{aligned}
{\check{\delta}}_{\b,n+1}:=\max\{\delta\leq &\tilde{\delta}_{\b,n} \ : \ \frac{\Phi(\ell,\b,\delta)}{\|\a\cdot\ell\|} g(c,n,\ell)  \leq \max_{\bar{\b}\in\Bbb{S}^{d-1}}\Phi(\ell,\bar{\b},\tilde{\delta}_{\bar{\b},n}) \\
&\quad\forall\ 0< |\ell|\leq 2^{n+1}\}.
\end{aligned}
\]

From the fact that $\tilde{\delta}_{\b,n}\leq \delta_{\b,n}$ for all $\b\in\Bbb{S}^1$ and using also Lemma \ref{facciamolo} we have 
$\bar{\delta}_{\b,n+1}\leq \delta_{\b,n+1}$ for all $\b\in\Bbb{S}^1$, and hence it suffices to show that $\tilde{\delta}_{\b,n+1}\geq\bar{\delta}_{\b,n+1}$. 
Now let us start by that if $\b\in\Bbb{S}^1$ and 
\[
\theta(\b,\tilde{\nu}(l+1)),\theta(\b,-\tilde{\nu}(l+1))\geq {C_n}{t_k^{-\frac{\bar{\theta}}{2}}}, 
\]
then for $\ell\in\Bbb{Z}^2$ such that $\ell\parallel \tilde{\nu}(l+1)$ and $|\ell|\leq 2^{n+1}$ one has
\beq{coseno}
\frac{\Phi(\ell,\b,\tilde{\delta}_{\b,n})}{\|\a\cdot\ell\|}\leq\Phi(\ell,\b_{l+1},\tilde{\delta}_{\b_{l+1},n}),
\qquad\qquad \b_{l+1}:=\frac{\tilde{\nu}(l+1)}{|\tilde{\nu}(l+1)|_2}
\eeq
Indeed one has
\[
\Phi(\ell,\b,\tilde{\delta}_{\b,n})=e^{2\pi\tilde{\delta}_{\b,n} |\b\cdot\ell|}=e^{2\pi\tilde{\delta}_{\b,n} |\ell|_2\cos\tilde\theta}
\]
where $\tilde\theta=\theta(\b,\tilde{\nu}(l+1))$. Then, from the fact that
\[
\theta(\b,\tilde{\nu}(l)),\theta(\b,-\tilde{\nu}(l))\geq {C}{t_k^{-\frac{\bar{\theta}}{2}}},
\]
we obtain
\[
\cos\tilde\theta\leq 1-C\tilde\theta^2\leq 1-{C}{-t_k^{\bar{\theta}}},
\]
from which we deduce
\beq{cos1}
\Phi(\ell,\b,\bar{\delta}_{\b,n})\leq e^{2\pi\bar{\delta}_{\b,n} |\ell|_2(1-{C}{t_k^{-\bar{\theta}}})},
\eeq
so that by \eqref{ncora}, \eqref{ncora1} we obtain \eqref{coseno}.

On the other hand, from the fact that $\ell\parallel \tilde{\nu}(l+1)$ we have
\beq{cos2}
\Phi(\ell,\b_{l+1},\bar{\delta}_{\b_{l+1},n})=e^{2\pi\bar{\delta}_{\b_{l+1},n}|\ell|_2}.
\eeq
Moreover, by \eqref{ncora}, \eqref{ncora1}, and Lemmata \ref{log} and \ref{x} we get, for $\b\in\Bbb{S}^1$ as in case 1,
 %such that there exists $t_k\leq h\leq l+1$ such that $\theta(\b,\tilde{\nu(h)})<\frac{C}{t_k^{\frac{\bar{\theta}}{2}}}$:
\[
\begin{aligned}
\max&\{\delta\leq \tilde{\delta}_{\b,n} \ :  \ \frac{\Phi(\ell,\b,\delta)}{\|\a\cdot\ell\|}\leq \Phi(\ell,\b_{l+1},\tilde{\delta}_{\b_{l+1},n}),\ |\ell|\leq 2^{n+1}\} \\
&\leq
\max\{\delta\leq \tilde{\delta}_{\b_{l+1},n} \ : \ \frac{\Phi(\ell,\b,\delta)}{\|\a\cdot\ell\|}\leq \Phi(\ell,\b_{l+1},\tilde{\delta}_{\b_{l+1},n}),\ |\ell|\leq 2^{n+1}\}-c_{n+1}\\
&\leq
\max\{\delta\leq \tilde{\delta}_{\b_{l+1},n} \ :  \ \frac{\Phi(\ell,\b_{l+1},\delta)}{\|\a\cdot\ell\|}\leq \Phi(\ell,\b_{l+1},\tilde{\delta}_{\b_{l+1},n}),\\
&\qquad\qquad \ell\parallel \tilde{\nu}(l+1), |\ell|\leq 2^{n+1}\}-c_{n+1}\\
&\leq \tilde{\delta}_{\b_{l+1},n}-C\Big(\max_{n-(\log k_n)^{\t}\leq p \leq n}\frac{p_n}{2^{p-n(t_k)} t_k\log t_k}-c_{n+1}\Big)\\
&\leq \tilde{\delta}_{\b_{l+1},n}-C\Big(\frac{(\log t_k)^{\t}p_n}{t_k 2^{n-n(t_k)}}-c_{n+1}\Big),
\end{aligned}
\]
which in turn implies that, for $\b\in\Bbb{S}^1$ as in case 1, one has
%so that we get for $\b\in\Bbb{S}^1$ such that there exists $t_k\leq h\leq l+1$ such that $\theta(\b,\tilde{\nu(h)})<\frac{C_{n}}{t_k^{\frac{\bar{\theta}}{2}}}$:
\[
\check{{\delta}}_{\b,n+1}\geq \delta-C\sum_{m=n(t_k)}^{ n+1}\big(\frac{(\log t_k)^{\t}p_m}{2^{m-n(t_k)}t_k}- c_m\big).
\]

Finally, for $\b\in\Bbb{S}^1$ as in case 2, using the fact that for $|\ell|\leq 2^{n+1}$ one has
%from the fact that for $\b\in\Bbb{S}^1$ such that there does not exists $t_k\leq h\leq l+1$ such that $\theta(\b,\tilde{\nu(h)})<\frac{C_{n}}{t_k^{\frac{\bar{\theta}}{2}}}$ and for $|\ell|\leq 2^{n+1}$:
\[
\frac{\Phi(\ell,\b,\delta_{\b,n}-c_{n+1})}{\|\a\cdot\ell\|}\leq \Phi(\ell,\b_{l+1},\bar{\delta}_{\b_{l+1},n}),
\]
we obtain
\[
\check{{\delta}}_{n+1,\b}\geq \delta- C\sum_{m=n(t_k)}^{ n+1}c_m.
\]
%%In particular, for $\b\in\Bbb{S}^1$ such that there exists $t_k\leq h\leq l+1$ such that $\theta(\b,\tilde{\nu(h)})<\frac{C_{n+1}}{t_k^{\frac{\bar{\theta}}{2}}}$:
%%\[
%%\bar{\delta}_{n+1,\b}\geq \delta-\sum_{n(t_k)\leq m\leq n+1}\frac{C(\log t_k)^{\t}p_m}{2^{m-n(t_k)}t_k}- C\sum_{n(t_k)\leq m\leq n+1}c_m.
%%\]
%%and for $\b\in\Bbb{S}^1$ such that there does not exists $t_k\leq h\leq l+1$ such that $\theta(\b,\tilde{\nu(h)})<\frac{C_{n+1}}{t_k^{\frac{\bar{\theta}}{2}}}$:
%%\[
%%\bar{\delta}_{n+1,\b}\geq \delta- C\sum_{n(t_k)\leq m\leq n+1}c_m.
%%\]
%%
In both cases the assertion follows.
\qed

\section{Construction of an explicit example}

We are now in the position to construct an explicit example of a vector $\a\in\Bbb{R}^2$ satisfying the Hypotheses of
Proposition \ref{criterio}. Recall Definition \ref{ricordatela}.

\begin{definition}\label{phi}
For $n\in\Bbb{N}$, consider any $l\in\Bbb{N}$ such that $t_n<l\leq t_{n+1}$ and set
    $$\Phi(l):=\frac{2t_n^{\frac{\bar{\theta}}{2}}}{l-t_n}-1.$$
\end{definition}
The r\^ole of the $\Phi(l)$ above is to twist the worst direction when we change the scale.

We define $\a=(\a_1,\a_2)\in\Bbb{R}^2$ as the limit of a sequence $(\a_{1,n},\a_{2,n})$ constructed iteratively in the following way.
Set $p_0=\bar{p}_0:=1$, $q_0:=100$, $\bar{q}_0:=101$ and $\bar{q}_{-1}:=1$.
Then define
\beq{a}
\left\{
\begin{array}{l}
\dst \a_{1,0}:=\frac{p_0}{q_0} ,\\  \ \\
\dst \a_{2,0}:=\frac{\bar{p}_0}{\bar{q}_0}.
\end{array}\right.
\eeq

Next, consider $l\in\Bbb{N}$ and suppose that we recursively defined $\a_{1,2l},\a_{2,2l}$ as
\beq{b}
\left\{
\begin{array}{l}
\dst \a_{1,2l}:=\frac{p_{2l}}{q_{2l}} ,\\  \ \\
\dst \a_{2,2l}:=\frac{\bar{p}_{2l}}{\bar{q}_{2l}},
\end{array}
\right.
\eeq
for some $p_0,\ldots,p_{2l}$, $q_0,\ldots,q_{2l}$, $\bar{p}_0,\ldots,\bar{p}_{2l}$, $\bar{q}_0,\ldots,\bar{q}_{2l}$
satisfying
\beq{iter}
\dst \bar{q}_{2l-1}^2q_{2l}\leq \bar{q}_{2l}\leq q_{2l} \bar{q}_{2l-1}^4, \qquad{\rm{GCD}}(q_{2l},\bar{q}_{2l})\leq \bar{q}_{2l-1}.
\eeq
Note that \eqref{iter} is satisfied for $l=0$.
Then, setting
\beq{cl}
c(l):=\left\lfloor\frac{\bar{q}_{2l}}{q_{2l}}\Phi(l)\right\rfloor,
\eeq
i.e.~the integer part of ${\bar{q}_{2l}}\Phi(l)/q_{2l}$,
and (recall \eqref{normezeta})
\beq{nu}
\tilde{\nu}(l):=(q_{2l}c(l),-\bar{q}_{2l}), \quad \nu(l):=\frac{\tilde{\nu}(l)}{{\rm{GCD}}(q_{2l}c(l),-\bar{q}_{2l})}, \quad \tilde{q}_{l}:=|\nu(l)|,
\eeq
 and denoting
 \beq{a2l1}
a_{2l+1}:=\left\lfloor e^{\frac{\tilde{q}_{l}}{t_n^2}}\right\rfloor,
\eeq
with $n$ such that $t_n<l\leq t_{n+1}$, we define
\beq{c}
\left\{
\begin{array}{l}
\dst \a_{1,2l+1}:=\a_{1,2l}+\frac{1}{a_{2l+1}q_{2l}c(l)} ,\\  \ \\
\dst \a_{2,2l+1}:=\a_{2,2l}+\frac{1}{a_{2l+1}\bar{q}_{2l}}.
\end{array}
\right.
\eeq
%
%\beq{}
%\tilde{\nu}(l):=(q_{2l}c(l),-\bar{q}_{2l}), \quad \nu(l):=\frac{\tilde{\nu}(l)}{{\rm{MCD}}(q_{2l}c(l),-\bar{q}_{2l})}, \tilde{q}_{l}:=|\nu(l)|,
%\eeq
%    where $|(l_1,l_2)|:=\max\{|l_1|,|l_2|\}$, and:
%\beq{}
%a_{2l+1}:=\left\lfloor e^{\frac{\tilde{q}_{l}}{t_n^2}}\right\rfloor,
%\eeq
%with $t_n<l\leq t_{n+1}$.

Note that by (\ref{iter}) one has $q_{2l}\leq \bar{q}_{2l}$ so that $c(l)\not=0$, which in turn implies that $\a_{1,2l+1}$ is well defined. 
Note also that the bigger $c(l)$ is, the closer is the direction of $\tilde{\nu}(l)$ to $(1,0)$. In particular the function $\Phi$ in Definition \ref{phi}
ensures that for different $l$ the direction of $\tilde{\nu}(l)$ changes.
Next we set
\[
\begin{aligned}
q_{2l+1}&:=a_{2l+1}q_{2l}^2c(l), \quad \bar{q}_{2l+1}:=a_{2l+1}\bar{q}_{2l}^2,\\
p_{2l+1}&:=q_{2l+1}p_{2l}+q_{2l}, \quad \bar{p}_{2l+1}:=\bar{q}_{2l+1}\bar{p}_{2l}+\bar{q}_{2l},
\end{aligned}
\]
so that one has $\a_{1,2l+1}=\frac{p_{2l+1}}{q_{2l+1}}$, $\a_{2,2l+1}=\frac{\bar{p}_{2l+1}}{\bar{q}_{2l+1}}$.
Finally, we define
\beq{d}
\left\{
\begin{array}{l}
\dst \a_{1,2l+2}:=\a_{1,2l+1}+\frac{1}{a_{2l+2}} ,\\  \ \\
\dst \a_{2,2l+2}:=\a_{2,2l+1}+\frac{1}{a_{2l+2}\bar{q}_{2l+1}q_{2l+1}+1},
\end{array}
\right.
\eeq
with 
\beq{a2l2}
\begin{aligned}
a_{2l+2}:=&\left\lfloor e^{\frac{\tilde{q}_{l}}{t_n\log t_n}} \right\rfloor,\qquad
q_{2l+2}:=a_{2l+2}q_{2l+1}, \\
 &\quad\bar{q}_{2l+2}:=(a_{2l+2}q_{2l+1}\bar{q}_{2l+1}+1)\bar{q}_{2l+1}.
 \end{aligned}
\eeq

Note that if $p>1$ is such that $p$ divides both $q_{2l+2}$ and $\bar{q}_{2l+2}$, then $p$ divides $\bar{q}_{2l+1}$ as well, because $p$ cannot 
divide $a_{2l+2}q_{2l+1}\bar{q}_{2l+1}+1$. 
This implies ${\rm{GCD}}(\bar{q}_{2l+2},q_{2l+2})\leq \bar{q}_{2l+1}$. Moreover, by \eqref{a2l2} and the definition of $a_{2l+1}$, $a_{2l+2}$, 
it follows that $\bar{q}_{2l+1}^2q_{2l+2}\leq \bar{q}_{2l+2}\leq q_{2l+2} \bar{q}_{2l+1}^4$. In other words \eqref{iter} is satisfied also up to $l+1$ and hence 
we can iterate the construction above. Finally we define
\beq{elui!}
\a=(\a_1,\a_2):=\lim_{l\to\infty} (\a_{1,2l},\a_{2,2l}).
\eeq

%In order to prove Theorem \ref{super} we need to prove that both components $\a_1,\a_2$ of $\a$ above satisfy \eqref{bry} with $d=1$ but $\a$ is not a Bryuno vector.
Frist note that the sizes of $\tilde{q}_l$ tand $\tilde{\nu}(l)$ are comparable because, since $q_{l}$ and $\bar{q}_{l}$ grow fast, then
 ${\rm{GCD}}(q_{2l},\bar{q}_{2l})$ is quite small; precisely we have the following result.
 
\lem{primastimetta}
For all $l\in\Bbb{N}$, one has $\tilde{q}_l\geq\frac{|\tilde{\nu}(l)|}{2l\bar{q}_{2l-1}^5}$.
\elem
\proof
By (\ref{iter}) one has ${\rm{GCD}}(q_{2l},\bar{q}_{2l})\leq \bar{q}_{2l-1}$. Moreover, by Definition \ref{phi} and \eqref{cl} we have 
$c(l)\leq \Phi(l) \bar{q}_{2l-1}^4\leq 2\sqrt{t_n} \bar{q}_{2l-1}^4\leq 2l \bar{q}_{2l-1}^4$, since 
$t_n<l$. But then, 
$$
{\rm{GCD}}(q_{2l}c(l),-\bar{q}_{2l})\leq c(l){\rm{GCD}}(q_{2l},\bar{q}_{2l})\leq \bar{q}_{2l-1}c(l)\leq 2l\bar{q}_{2l-1}^5,
$$
so the assertion follows.
 \qed

Note that $q_l$ and $p_l$ could have some common factor. However their greatest common divisor is not ``too big'' (same applies to $\bar{q}_l$ and $\bar{p}_l$). Precisely 
we have the following result.

\lem{urca}
For all $l\in\Bbb{N}$ one has ${\rm{GCD}}(q_l,p_l)\leq q_{l-1}$ and ${\rm{GCD}}(\bar{q}_l,\bar{p}_l)\leq \bar{q}_{l-1}$.
\elem

\proof
Let $p\in\Bbb{N}$ such that $p$ divides both $p_l$ and $q_l$.Then $p$ divides $p_{l}-q_{l}p_{l-1}$ i.e. it divides $q_{l-1}$. The proof for $\bar{p}_l,\bar{q}_l$ is the same. 
\qed

\lem{3}
For all $l\in\Bbb{N}$ one has
\[
|\a_{i,2l+1}-\a_i|\leq \frac{2}{a_{2l+2}},\qquad i=1,2.
\]
\elem

\proof
By definition of $\a_1,\a_2$ we have 
\[
|\a_{i,2l+1}-\a_i| \leq \sum_{j=2l+2}^{\infty}\su{a_{j}},\qquad i=1,2.
\]
Then, the result follows from the fact that for any $j\ge1$ one has $a_{j+1}>a_{j}^2$.
\qed

We now want to give upper and lower bounds on %the loss of regularity at scale $\tilde{q}_l$, %. In fact, once we 
 $\Omega_{\a}(\tilde{q}_l)$ (recall Definition \ref{bryvec}), 
since the loss of regularity at scale $\tilde{q}_l$ 
is of size $\sim\frac{1}{\tilde{q}_l}\log(\su{\Omega_{\a}(\tilde{q}_l)})$.

\lem{stoqua}
For all $l\in\Bbb{N}$ one has
\beq{uplow}
\su{2(a_{2l+2}\bar{q}_{2l+1}q_{2l+1}+1)}\leq \Omega_{\a}(\tilde{q}_l)\leq \frac{8\bar{q}_{2l+1}^5l}{a_{2l+2}}.
\eeq
\elem

\proof
First of all notice that, given any $\ell=(\ell_1,\ell_2)\in\Bbb{Z}^2$, one has
\[
\begin{aligned}
 \| \a \cdot \ell \|&\leq \|\a_{1,2l+1}\ell_1+\a_{2,2l+1}\ell_2 \| \\
 &\qquad+
 \Bigl| (\a_1-\a_{1,2l+1}) \ell_1+(\a_2-\a_{2,2l+1}) \ell_2 \Bigr|,
 \end{aligned}
 \]
so that using Lemma \ref{3}, if $|\ell |\leq \tilde{q}_l$ one has
\[
 \Bigl| (\a_1-\a_{1,2l+1}) \ell_1+(\a_2-\a_{2,2l+1}) \ell_2 \Bigr|\leq \frac{4(|\ell_1|+|\ell_2|)}{a_{2l+2}}\leq\frac{4\tilde{q}_l}{a_{2l+2}},
\]
and hence one has
\beq{lasagna}
\begin{aligned}
 \| \a \cdot \ell \|&\leq \frac{4\tilde{q}_l}{a_{2l+2}}+\|(\a_{1,2l+1},\a_{2,2l+1})\cdot\ell \|.
 \end{aligned}
\eeq
On the other hand one has
\beq{minimo}
\min_{|\ell|\le \tilde{q}_l} \|(\a_{1,2l+1},\a_{2,2l+1})\cdot\ell \| = \|\a_{1,2l+1}\hat{\ell}_1+\a_{2,2l+1}\hat{\ell}_2\|=0
\eeq
where (recall \eqref{nu}) we denoted
\[
\hat{\ell}_1:= \frac{q_{2l}c(l)}{{\rm{GCD}}(q_{2l}c(l),-\bar{q}_{2l})},
\qquad
\hat{\ell}_2:= - \frac{\bar{q}_{2l}}{{\rm{GCD}}(q_{2l}c(l),-\bar{q}_{2l})}.
\]
In turn \eqref{lasagna} and \eqref{minimo} imply that
$$
\Omega_{\a}(\tilde{q}_{l})=\min_{|\ell|\le \tilde{q}_l}\| \a \cdot \ell \|\leq \frac{4\tilde{q}_l}{a_{2l+2}}\leq\frac{4c(l)\bar{q}_{2l}}{a_{2l+2}}\leq \frac{8l\bar{q}_{2l+1}^5}{a_{2l+2}},
$$
where the last inequality follows from the fact that $c(l)\leq 2l\bar{q}_{2l+1}^4$ and $\bar{q}_{2l}\leq \bar{q}_{2l+1}$.

\giu
\nl
In order to prove the lower bound in \eqref{uplow} we consider any $\ell=(\ell_1,\ell_2)$ with $|\ell|\le\tilde{q}_{l}$ and study separately two cases.

\begin{itemize}

\item
{\it{Case 1: the vector $\ell$ is not parallel to $\nu(l)$.}} In this case we write
\[
\begin{aligned}
\|\a \cdot\ell \| &\geq \Big| \|\a_{1,2l+1}\ell_1+\a_{2,2l+1}\ell_2\|-|(\a_1-\a_{1,2l+1})\ell_1|  \\
 &\qquad-|(\a_2-\a_{2,2l+1})\ell_2| \Big|
\end{aligned}
\]
where, by the definition of $\a_{1,2l+1}$ and $\a_{2,2l+1}$ in \eqref{c},
\[
\|\a_{1,2l+1}\ell_1+\a_{2,2l+1}\ell_2\|=\left\|\left(\frac{p_{2l}}{q_{2l}}\ell_1+\frac{\bar{p}_{2l}}{\bar{q}_{2l}}\ell_2\right)+\frac{\ell_1}{c(l)q_{2l}a_{2l+1}}+\frac{\ell_2}{\bar{q}_{2l}a_{2l+1}}\right\|.
\]

\begin{enumerate}

\item
If $\frac{p_{2l}}{q_{2l}}\ell_1+\frac{\bar{p}_{2l}}{\bar{q}_{2l}}\ell_2$ is an integer we obtain
\beq{prima}
\begin{aligned}
\|\a_{1,2l+1}\ell_1+\a_{2,2l+1}\ell_2\|&= \left\|\frac{\ell_1}{c(l)q_{2l}a_{2l+1}}+\frac{\ell_2}{\bar{q}_{2l}a_{2l+1}}\right\|\\
&\geq\su{a_{2l+1}q_{2l}\bar{q}_{2l}c(l)},
\end{aligned}
\eeq
where the last inequality follows from the fact that $\ell$ is not parallel to $\nu(l)$.

\item 
If $\frac{p_{2l}}{q_{2l}}\ell_1+\frac{\bar{p}_{2l}}{\bar{q}_{2l}}\ell_2$ is not an integer we write
\beq{staqua}
\begin{aligned}
\|\a_{1,2l+1}\ell_1+&\a_{2,2l+1}\ell_2\| \geq 
\Big|\|\a_{1,2l}\ell_1+\a_{2,2l}\ell_2\| \\ 
&-|(\a_{1,2l+1}-\a_{1,2l})\ell_1|-|(\a_{2,2l+1}-\a_{2,2l})\ell_2|\Big|,
\end{aligned}
\eeq
where, again by Lemma \ref{3}, we bound
\beq{P1}
|(\a_{1,2l+1}-\a_{1,2l})\ell_1|+|(\a_{2,2l+1}-\a_{2,2l})\ell_2|\leq\frac{4\tilde{q}_l}{a_{2l+1}q_{2l}} \le \frac{8l\bar{q}_{2l}\bar{q}_{2l-1}^4}{a_{2l+1}}
\eeq
while, since $\frac{p_{2l}}{q_{2l}}\ell_1+\frac{\bar{p}_{2l}}{\bar{q}_{2l}}\ell_2$ is not an integer, we have
\beq{P2}
\|\a_{1,2l}\ell_1+\a_{2,2l}\ell_2\|=\left\|\frac{p_{2l}}{q_{2l}}\ell_1+\frac{\bar{p}_{2l}}{\bar{q}_{2l}}\ell_2\right\|\geq\frac{1}{q_{2l}\bar{q}_{2l}},
\eeq
so that using (\ref{P1}) and (\ref{P2}) into \eqref{staqua} we deduce
\[
\|\a_{1,2l+1}\ell_1+\a_{2,2l+1}\ell_2\|\geq\frac{1}{q_{2l}\bar{q}_{2l}}-\frac{8l\bar{q}_{2l}\bar{q}_{2l-1}^4}{a_{2l+1}} \geq\su{a_{2l+1}q_{2l}\bar{q}_{2l}c(l)}.
\]

\end{enumerate}

Overall in both cases we deduce that, if $\ell$ is not parallel to $\nu(l)$, then
\beq{questo}
\begin{aligned}
\|\a \cdot\ell \| &\geq \Big|\|\a_{1,2l+1}\ell_1+\a_{2,2l+1}\ell_2\|-|(\a_1-\a_{1,2l+1})\ell_1| \\
 &\qquad-|(\a_2-\a_{2,2l+1})\ell_2| \Big| \\
&\geq \su{a_{2l+1}q_{2l}\bar{q}_{2l}c(l)}-\frac{4\tilde{q}_l}{a_{2l+2}q_{2l}}\geq \su{2a_{2l+1}q_{2l}\bar{q}_{2l}c(l)}\geq \su{a_{2l+2}}.
\end{aligned}
\eeq

\item
{\it{Case 2: the vector $\ell$ is parallel to $\nu(l)$.}} In this case we have
\[
\begin{aligned}
\|\a\cdot\ell\|&=\|\a_{1,2l+2}\ell_1+ \a_{2,2l+2} \ell_2\|\\
&\geq \frac{(a_{2l+2}\bar{q}_{2l+1}q_{2l+1}+1)\ell_1+a_{2l+2} \ell_2 }{a_{2l+2}(a_{2l+2}\bar{q}_{2l+1}q_{2l+1}+1)}-\frac{8l\bar{q}_{2l}\bar{q}_{2l-1}^4}{a_{2l+3}}.
\end{aligned}
\]
But then, since ${\rm{GCD}}(a_{2l+2}\bar{q}_{2l+1}q_{2l+1}+1,a_{2l+2})=1$, we have
\[
\frac{(a_{2l+2}\bar{q}_{2l+1}q_{2l+1}+1)\ell_1+ a_{2l+2}\ell_2}{a_{2l+2}(a_{2l+2}\bar{q}_{2l+1}q_{2l+1}+1)}\geq\frac{1}{a_{2l+2}\bar{q}_{2l+1}q_{2l+1}+1},
\]
so that, by the fact that $|\ell|\leq 2l\bar{q}_{2l}\bar{q}_{2l-1}^4$, we obtain
\[
\|\a\cdot\ell\|\geq \frac{1}{2(a_{2l+2}\bar{q}_{2l+1}q_{2l+1}+1)}.
\]
\end{itemize}

\noindent
In both cases the assertion follows.
\qed

\rem{storem}
By Lemma \ref{stoqua} above we deduce that the loss of regularity at scale $\tilde{q}_l$ is of order $\su{t_n\log t_n}$ with $n$ such that $t_n<l\leq t_{n+1}$.
Indeed one has
\[
\frac{1}{\tilde{q}_l} \log \Omega_{\a}^{-1}(\tilde{q}_l)\leq\frac{\log 2(a_{2l+2}\bar{q}_{2l+1}q_{2l+1}+1)}{\tilde{q}_l}\leq\frac{4\log a_{2l+2}}{\tilde{q}_l}\leq \frac{4}{t_n\log t_n},
\]
where the last inequality follows from the definition of $a_{2l+2}$ in \eqref{a2l2}.
On the other hand one has
\beq{thebound}
\begin{aligned}
\frac{1}{\tilde{q}_l} \log \Omega_{\a}^{-1}(\tilde{q}_l) &\geq \frac{1}{\tilde{q}_l} \log \frac{a_{2l+2}}{8l \bar{q}_{2l+1}^5}=\frac{1}{\tilde{q_l}}( \log a_{2l+2}-{\log 8l\bar{q}_{2l+1}^5 })\\
&\geq\frac{1}{2\tilde{q}_l}  \log a_{2l+2}\geq\frac{1}{2t_n\log t_n}.
\end{aligned}
\eeq
\erem

\rem{remarkino}
The bound \eqref{thebound} in Remark \ref{storem} tells us that
\[
\Omega_\a(\tilde{q}_l) \le e^{-\frac{\tilde{q}_l}{2t_n\log t_n}}.
\]
On the other hand the bound \eqref{questo} implies that
\[
\min_{\substack{ \ell\in\Bbb{Z}^2 \\ \ell \nparallel \nu(l) \\ |\ell|\le \tilde{q}_l }} \| \a \cdot\ell \| \ge \frac{1}{2a_{2l+1}q_{2l}\bar{q}_{2l}c(l)}.
\]
Thus $\Omega_\a(\tilde{q}_l)$ is attained for some $\ell\parallel \nu(l)$. On the other hand by construction one has ${\rm GCD}(\nu(l)_1, \nu(l)_2)=1$ (recall \eqref{nu}), 
and hence $\Omega_\a(\tilde{q}_l)$ is attained either at $\nu(l)$ or at $-\nu(l)$. 
\erem

\lem{C1}
Let $l\in\Bbb{N}$ and consider any $\ell\in\Bbb{Z}^2$ such that $\tilde{q}_l\leq |\ell|_1<\sqrt{\tilde{q}_{l+1}}$, and $\ell$ is not parallel to $\nu(l)$. Then one has
\[
\frac{1}{|\ell|_1} \log\left(\su{\|\a\cdot\ell \|}\right) \leq \frac{3\tilde{q}_l}{t_n^2|\ell|_1},
\]
with $n$ such that $t_n<l\leq t_{n+1}$.
\elem

\proof
Since $\tilde{q}_l\leq |\ell|_1<\sqrt{\tilde{q}_{l+1}}$ and $\ell$ is not parallel to $\nu(l)$, 
we use Lemma \ref{3} once more to obtain
\[
\|\a\cdot\ell\|\geq \|\a_{1,2l+1} \ell_1+\a_{2,2l+1} \ell_2\|-\frac{4\sqrt{\tilde{q}_{l+1}}}{a_{2l+2}}\geq \frac{1}{2a_{2l+1}q_{2l}\bar{q}_{2l}}\geq\frac{1}{a_{2l+1}^3}.
\]
But then by the definition of $a_{2l+1}$ in \eqref{a2l1} one has
\[
\frac{1}{|\ell|_1} \log\left(\su{\| \a\cdot\ell \|}\right)\leq \frac{3\log a_{2l+1}}{\tilde{q}_l}\frac{\tilde{q}_l}{|\ell|_1}\leq \frac{3\tilde{q}_l}{t_n^2|\ell|_1},
\]
so that the assertion follows.
\qed

\lem{sottoradice}
Let $\ell=m\nu(l)\in\Bbb{Z}^2$, with $m$ such that $\tilde{q}_l\le|m \tilde{q}_l|<\sqrt{\tilde{q}_{l+1}}$. Then one has
\elem
\[
 \frac{1}{4|m| l\log l}\leq\frac{1}{|\ell|_1} \log\left(\su{\|\a\cdot\ell \|}\right)\leq \frac{4}{|m| l\log l}.
\]

\proof
First observe that 
\[
\| m \nu(l) \cdot \a \| = |m| \| \nu(l)\cdot \a\|
\]
for $|m|<\sqrt{\tilde{q}_{l+1}}<\Omega_\a(\tilde{q}_l)^{-1}$. Then combining Lemma \ref{stoqua} with Remark \ref{remarkino}
\[
\frac{1}{|m|\tilde{q}_l} \log\left(\frac{a_{2l+2}}{8\bar{q}_{2l+1}^5 l}\right) \leq\frac{1}{\tilde{q}_l} \log\frac{1}{\|m\nu(l)\cdot \a\|}\leq \frac{1}{|m|\tilde{q}_l} \log\left(\frac{a_{2l+2}}{\bar{q}_{2l}l}\right),
\]
from which the result follows directly by the definition of $a_{2l+2}$ in \eqref{a2l2}.
\qed

\rem{unbel}
Note that Lemmata \ref{C1} and \ref{sottoradice} provide an upper bound for all $\ell\in\Bbb{Z}^2$ with $\tilde{q}_l\leq |\ell|_1<\sqrt{\tilde{q}_{l+1}}$.
In particular if $\ell \parallel \nu(l)$, Lemma \ref{sottoradice} provides also a lower bound.
\erem

\lem{C2}
Let $l\in\Bbb{N}$ and consider any $\ell \in\Bbb{Z}^2$ such that $\sqrt{\tilde{q}_{l+1}}\leq |\ell|_1<\tilde{q}_{l+1}$.  Then one has
\[
\frac{1}{|\ell|_1} \log\left(\su{\|\a\cdot\ell\|}\right)  \leq \frac{3\sqrt{\tilde{q}_{l+1}}}{t_n^2|\ell|_1},
\]
with $n$ such that $t_n<l\leq t_{n+1}$.
\elem

\proof
As in the proof of Lemma \ref{C1} we use Lemma \ref{3} to deduce that, if $\ell$ is not parallel to $\nu(l)$,
\[
\|\a \cdot \ell\| \geq \|\a_{1,2l+2} \ell_1+\a_{2,2l+2} \ell_2\|-\frac{4\tilde{q}_{l+1}}{a_{2l+3}}\geq \frac{1}{4a_{2l+2}q_{2l}\bar{q}_{2l}} \geq \frac{1}{a_{2l+2}^3}.
\]
On the other hand if $\ell$ is parallel to $\nu(l)$, we write $\ell=m\nu(l)$ and, noting that 
\[
{\rm GCD}({a_{2l+2} ,  a_{2l+2}\bar{q}_{2l+1}q_{2l+1}+1 })=1, \quad\mbox{ and }\quad
a_{2l+2}( a_{2l+2}\bar{q}_{2l+1}q_{2l+1}+1) > \tilde{q}_{l+1},
\]
we deduce
\begin{equation}\label{cremina}
\frac{1}{a_{2l+2}}\ell_1 + \frac{1}{a_{2l+2}\bar{q}_{2l+1}q_{2l+1}+1}\ell_2  \notin \Bbb{Z},
\end{equation}
so that we can use once more Lemma \ref{3} to obtain, for $|\ell|\leq\tilde{q}_{l+1}$,
\beq{above}
\begin{aligned}
\|\a \cdot\ell \|&\geq \|\a_{1,2l+2} \ell_1+\a_{2,2l+2} \ell_2\|-\frac{4\tilde{q}_{l+1}}{a_{2l+3}} \\
&= 
\| \frac{1}{a_{2l+2}}\ell_1 + \frac{1}{a_{2l+2}\bar{q}_{2l+1}q_{2l+1}+1}\ell_2 \| -\frac{4\tilde{q}_{l+1}}{a_{2l+3}} \\
&\stackrel{\eqref{cremina}}{\geq} \frac{1}{ a_{2l+2}( a_{2l+2}\bar{q}_{2l+1}q_{2l+1}+1)} -\frac{4\tilde{q}_{l+1}}{a_{2l+3}}  \geq \frac{1}{a_{2l+2}^3}.
\end{aligned}
\eeq

But then by the definition of $a_{2l+2}$ in \eqref{a2l2} we have
\[
\frac{1}{|\ell|_1} \log\left(\su{\|\a\cdot\ell \|}\right)  \leq \frac{3\log a_{2l+2}}{\sqrt{\tilde{q}_{l+1}}}\frac{\sqrt{\tilde{q}_{l+1}}}{|\ell|_1}\leq \frac{3\sqrt{\tilde{q}_{l+1}}}{t_n^2|\ell|_1},
\]
which proves the assertion.
\qed

\pro{nobruno}
The vector $\a=(\a_1,\a_2)$ defined in \eqref{elui!} is not a Bryuno vector.
\epro

\proof
One has
\[
\sum_{k\geq1}\frac{1}{2^k} \log \su{\Omega_{\a}(2^k)}\geq \sum_{l\geq1} \frac{1}{2\tilde{q}_l} \log\su{\Omega_{\a}(\tilde{q}_l)}\geq \sum_{l\ge1}\su{8l\log l}=+\infty,
\]
see also Remark \ref{remarkino}, so that $\a$ is not a Bryuno vector.
\qed

In order to conclude the proof of Theorem \ref{super} we need to prove that $\a$ is weak-Bryuno vector.

\lem{prossimo}
Fix $\b\in\Bbb{S}^1$ and suppose that there is $l$ large enough such that 
\[
\min\{|\b - \frac{\nu(l)}{|\nu(l)|_2}| , |\b + \frac{\nu(l)}{|\nu(l)|_2}| \} < {c}{l^{-\frac{\bar{\theta}}{2}}},
\]
for some $c>0$.
Let $n$ such that $t_n<l\le t_{n+1}$. Then for any $j=l+1,\ldots, l + \lfloor \frac{t_n}{4}\rfloor$ one has
\[
\min\{|\b - \frac{\nu(j)}{|\nu(j)|_2}| , |\b + \frac{\nu(j)}{|\nu(j)|_2}| \} \ge {c}{l^{-\frac{\bar{\theta}}{2}}}.
\]
\elem

\proof
First of all note that, by construction,
\[
\tilde{\nu}(l)= \bar{q}_{2l}(\Phi(l),-1) - (q_{2l}\left\{ \frac{\bar{q}_{2l} \Phi(l)}{q_{2l}} \right\},0)
\]
and hence
\[
\frac{\nu(l)}{|\nu(l)|_1} = \frac{1}{{\Phi(l)+1}} (\Phi(l),-1) + \e_l,\qquad |\e_l|\le\su{\bar{q}_{2l-1}^2}.
\]
But then,  for $j=l+1,\ldots, l + \lfloor \frac{t_n}{4}\rfloor$, one has, for some constant $c$,
\[
\begin{aligned}
\left| \frac{\nu(l)}{|\nu(l)|_1}  - \frac{\nu(j)}{|\nu(j)|_1}\right|   &\ge c\left | \frac{1}{{\Phi(l)+1}} - \su{{\Phi(j)+1}}\right|  \\
& \ge {c}{t_n^{-\frac{\bar{\theta}}{2}}} \ge {c}{l^{-\frac{\bar{\theta}}{2}}}.
\end{aligned}
\]
In particular, there exists $C>0$ such that for $j=l+1,\ldots, l + \lfloor \frac{t_n}{4}\rfloor$ one has
\[
\begin{aligned}
\left| \frac{\nu(l)}{|\nu(l)|_2}  - \frac{\nu(j)}{|\nu(j)|_2}\right|   &\ge C\left | \frac{1}{{\Phi(l)+1}} - \su{{\Phi(j)+1}}\right|  \\
& \ge {c}{t_n^{-\frac{\bar{\theta}}{2}}} \ge {c}{l^{-\frac{\bar{\theta}}{2}}}.
\end{aligned}
\]
and then the assertion follows by Remark \ref{storem}. 
\qed

\medskip

\noindent
{\bf{Proof (of Theorem \ref{super})}}
We will show that the vector $(\a_1,\a_2)$ satisfies the properties of Proposition \ref{criterio}. To show that property 1 holds, note that, by Lemma \eqref{primastimetta} and by \eqref{a2l1} we have:
\[
\tilde{q}_{l+1}\geq \frac{a_{2l+1}\tilde{q}_l}{2l \tilde{q}_{l-1}^5}\geq 4\tilde{q}_l.
\]
In particular, the first property holds.
The second property follows by Remark \ref{storem}.
The third property follows by Lemma \ref{prossimo}.
The fourth property follows by Lemma \ref{C1} and \ref{C2}.
The last property follows by Lemma \ref{sottoradice} and Lemma \ref{C2}.
\qed

%%%%%%%%%%%%%%%%%%%%%%%%%%%%%%%%%%%%%%%%%%%%%%%%%%%%%%%%%%%%%%%%%%%%%%%%%%%%%%%%%%%%%%%%%%%%

\zerarcounters

\section{The KAM scheme}\label{provakam}

In order to prove Theorem \ref{main} we take as fundamental domain for the torus the cube $[0,1]^d$, 
and we parametrize it with segments from the origin to the boundary of the cube.
Then, according to the scale in the KAM step, we will loose regularity not in a uniform way, but differently on each segment.
 In particular, the domain in which we will have a good control at each KAM step will not be some strip as usual. 
 
 As it is usual in KAM theory, one needs a to consider a ``thickened'' torus. So we give the following definition.
 
 \begin{definition}\label{param}
 Let $\xi:=\{\delta_\b\}_{\b\in\Bbb{S}^{d-1}}$ be such that $\inf_{\b\in \Bbb{S}^{d-1}}\delta_{\b}>0$, and set
 \[
 \Bbb{T}^d_{\xi}:= \bigcup_{\b\in\Bbb{S}^{d-1}} S_{\b,\delta_\b}
 \]
 where, given any $\b\in\Bbb{S}^{d-1}$ and any $\bar{\delta}>0$ we denoted
 \beq{segmento}
S_{\b,\bar{\delta}}:=
\{x+i\b y \;:\; x\in\Bbb{T}^d, y\in\Bbb{R},\,0\leq |y|<\bar{\delta}\}
 \eeq
 
 \end{definition}

We now introduce a suitable norm for analytic functions defined on   $\Bbb{T}^d_{\xi}$,

\begin{definition}\label{ananorm}
With the same notation as in Definition \ref{param}, we say that $f\in D^{\omega}(\Bbb{T}^d_{\xi})$ if $f\in D^{0}(\Bbb{T}^d)$ and
$f$ is analytic on $\bigcup_{\b\in\Bbb{S}^{d-1}}S_{\b,\delta_\b}$. For $\b\in\Bbb{S}^{d-1}$, $(x,y)\in S_{\b,\delta_{\b}}$ let $f_{\b}(x,y):=f(x+i\by)$.
We then define%% ho spostato il sup dentro%%
\beq{normaxi}
{
\| f \|_{ \xi }:=\sum_{\ell \in\Bbb{Z}^d}\sup_{z\in  \Bbb{T}^d_{\xi} }\left|\hat{f}(\ell)e^{2\pi i z\cdot\ell}\right|.
}
\eeq
Moreover, given any $\b\in\Bbb{S}^{d-1}$, for any $\bar{\delta}\leq \delta_\b$ we define 
\beq{stanormina} %%sup dentro%%
{
|f_{\b}|_{\bar{\delta}}:=\sum_{\ell\in\Bbb{Z}^d}\sup_{z\in S_{\b,\bar{\delta}}}\left|\hat{f}(\ell)e^{2\pi iz\cdot\ell}\right|.
}
\eeq

\end{definition}

\rem{deltaxi}
If $\xi:=\{\delta_\b\}_{\b\in\Bbb{S}^{d-1}}$ is such that $\delta_\b=\delta$ for all $\beta\in\Bbb{S}^{d-1}$ one has
\[
\|f\|_{\xi}=\sum_{\ell \in\Bbb{Z}^d}\sup_{\substack{z\in\Bbb{T}^d_\xi }}\left|\hat{f}(\ell )e^{2\pi iz\cdot\ell}\right| = \sum_{\ell\in\Bbb{Z}^d}|\hat{f}(\ell)|e^{2\pi \delta|\ell|_2},
\]
which is the standard analytic norm in \eqref{normana}.
\erem

\rem{supp}
Given any $\xi$ as in Definition \ref{param} and any $f\in D^{\omega}(\Bbb{T}^d_{\xi})$, one has
\[
\sup_{\b\in\Bbb{S}^{d-1}}|f_{\b}|_{\delta_{\b}}\leq\| f \|_{\xi}.
\]
\erem

{
The norm $\|\cdot\|_\xi$ in \eqref{normaxi} satisfies the algebra property as the next result shows.
\lem{prodotto}
Let $\xi=\{\delta_{\b}\}_{\b\in\Bbb{S}^{d-1}}$ as above. For any $f_1,f_2\in D^{\omega}(\Bbb{T}^d_\xi)$ one has
\[
\|f_1 f_2\|_{\xi}\leq \|f_1\|_{\xi}\|f_2\|_{\xi}.
\]

\elem
\proof
By definition we have
\[
\begin{aligned}
\|f_1 f_2\|_{\xi}&=\sum_{\ell\in\Bbb{Z}^d}\sup_{z\in S_{\b,\delta}}\left|\left(\sum_{\ell_1+\ell_2=\ell}\hat{f}_1(\ell_1)\hat{f}_2(\ell_2) 
e^{2\pi i z\cdot \ell}\right)\right| \\
&\leq \sum_{\ell\in\Bbb{Z}^d}\left(\sum_{\ell_1+\ell_2=l}\sup_{z\in S_{\b,\delta}}\left|\hat{f}_1(\ell_1)e^{2\pi i z\cdot \ell_1}\right| 
\sup_{z\in S_{\b,\delta}}\left|\hat{f}_2(\ell_2)e^{2\pi i z\cdot \ell_2}\right|\right)\leq \|f_1\|_{\xi}\|f_2\|_{\xi}
\end{aligned}
\]
so the assertion follows.
\qed
}

\begin{definition}\label{fourier}
    For any $N\in\Bbb{N}$ and any $g\in D^{\omega}(\Bbb{T}^d_\xi)$ we denote
    \[
    \begin{aligned}
    \T_{N}g(x)&:=\sum_{\substack{\ell\in\Bbb{Z}^d \\ |\ell |\leq N}}\hat{g}(\ell )e^{2\pi i x\cdot\ell},  \qquad\qquad
    \R_{N}g(x)&:=g(x)-\T_{N}g(x).
    \end{aligned}
    \]
\end{definition}

Then norm $\|\cdot\|_\xi$ in \eqref{normaxi} satisfies also the following.
{
\lem{derivate}
Let $N\in\Bbb{N}$, $\xi=\{\delta_{\b}\}_{\b\in\Bbb{S}^{d-1}}$ as in Definition \ref{param} 
and $f\in D^{\omega}(\Bbb{T}^d_\xi)$ such that $f(z)=\T_{N}f(z)$. Then, for $t_1,\dots t_d\in\Bbb{N}$ we have
\[
\left\|\frac{\partial^{t_1+\dots+t_d}f}{\partial_1^{t_1}\dots\partial_d^{t_d}}\right\|_{\xi}\leq (2\pi N)^{t_1+\dots+t_d}\|f\|_{\xi}.
\]
\elem
\proof
Since $f(z)=\T_{N}f(z)$ we have
\[
\begin{aligned}
\left\|\frac{\partial^{t_1+\dots+t_d}f}{\partial_1^{t_1}\dots\partial_d^{t_d}}\right\|_{\xi}&
=\sum_{\substack{l=(l_1,\dots,l_d)\in\Bbb{Z}^d \\ |l|\leq N}}\left|\sup_{z\in\Bbb{T}^d_{\xi}}|(2\pi l_1)^{t_1}|\dots |(2\pi l_d)^{t_d}|
|\hat{f}(l)|e^{2\pi i z\cdot l}\right| \\
&\leq (2\pi N)^{t_1+\dots+t_d}\|f\|_{\xi},
\end{aligned}
\]
so the assertion follows.
\qed
}
%%%
%%%\begin{definition}\label{mancounnome}
%%%  For any $n\in\Bbb{N}$ we define
%%%  \[
%%%    \theta_n:=\frac{|(\ell_n)_2|}{|\ell_n|_1},
%%%   \]
%%%    with $\ell_n=((\ell_n)_1,(\ell_n)_2)\in\Bbb{Z}^2$ as in Definition \ref{fernandovector}. Set also
%%%\[
%%%\g_n:=\inf\left\{\|\a\cdot y\| \, : \, y\in\Bbb{Z}^2, \, |\frac{y}{|y|_1}-\frac{\ell_n}{|\ell_n|_1}|\ge \sqrt{\frac{\log (\Omega_{\a}(2^n)^{-1})}{2^n}}\right\}
%%%\]
%%%    
%%%\end{definition}
%%%
%%%
%%%\lem{cista}
%%%One has
%%%\[
%%%\sum_{n\geq1}\su{\g_n}<\infty.
%%%\]
%%%\elem
%%%
%%%\proof
%%%
%%%\qed

\lem{Q}
Fix $\gotd>0$ and let $N$, $C$ and $c=\{c_n\}_{n\in\Bbb{N}}$ be such that \eqref{weakcond} holds. 
Then, there exists a sequence $\tilde{c}=\{\tilde{c}_n\}_{n\in\Bbb{N}}$ such that the following holds.
\begin{enumerate}
\item For every $C_1>0$ there exists $N_1\in\Bbb{N}$ such that, for $n\geq N_1$ one has $C_1 2^{nd}e^{-2^n\tilde{c}_n}<e^{-2^n c_n}$.
\item $\sum_{n\in\Bbb{N}}\tilde{c}_n<\infty$.
\end{enumerate}
Moreover, given any $\tilde{c}=\{\tilde{c}_n\}_{n\in\Bbb{N}}$ as above and for every $C_1>0$, there exists $N_1\in\Bbb{N}$ such that the following holds.
For $\b\in\Bbb{S}^{d-1}$, $0\leq n\leq N_1$ let $\delta_{\b,n}=1$.  For $n\geq N_1$ define recursively
\[
\tilde{\delta}_{\b,n+1}:=\max\{\delta\leq \delta_{\b,n}:\frac{\Phi(\ell,\b,\delta)}{\|\a\cdot\ell\|} g(c,n,\ell) \leq \max_{\bar{\b}\in\Bbb{S}^{d-1}}\Phi(\ell,\bar{\b},\delta_{\bar{\b},n})\quad\forall 2^n< |\ell|\leq 2^{n+1}\},
\]
and
\[
\delta_{\b,n+1}:=\min_{|\b-\bar{\b}|<e^{-{2^{n} \tilde{c}_n}}}\delta_{\bar{\b},n+1}.
\]
Then,
\[
\inf_{\b\in\Bbb{S}^{d-1}}\inf_{n\in\Bbb{N}}\delta_{\b,n}>0.
\]
\elem

\proof
Note that the first property is true if and only if $\tilde{c}_n\geq c_n+\frac{\log(C_12 ^{2nd})}{2^n}$. So, let:
\[
\tilde{c}_n:= c_n+\frac{\log(C_12 ^{2nd})}{2^n}.
\]
Then, by definition the first property holds, and from the fact that $(c_n)_{n\in\Bbb{N}}\in l^1$ we get that the second property holds.
Finally, from the fact that $\tilde{c}_n\geq c_n$ we get that by Lemma \ref{facciamolo} the sequence $(\bar{\delta}_{\b,n})$ defined in the same way as $(\delta)_{\b,n}$ (and with the same choice of $C,N$) but with $(c_n)_{n\in\Bbb{N}}$ instead of $(c_n)_{n\in\Bbb{N}}$ satisfies $\delta_{\b,n}\geq \bar{\delta}_{\b,n}$. From the fact that $\a$ is weak-Bryuno and by definition of $N,C, (c_n)_{n\in\Bbb{N}}$ we have that:
\[
\inf_{\b\in\Bbb{S}^{d-1}}\inf_{n\in\Bbb{N}}\delta_{\b,n}\geq \inf_{\b\in\Bbb{S}^{d-1}}\inf_{n\in\Bbb{N}}\bar{\delta}_{\b,n} >0.
\]

\qed

From now on, thanks to
Lemma \ref{Q} above, we can assume without loss of generality that the sequence $c=\{c_n\}_{n\in\Bbb{N}}$ satisfies the same properties 1-2 as the sequence $\{\tilde{c}_n\}_{n\in\Bbb{N}}$.
We now construct recursively a family of sequences parametrized by $\beta\in\Bbb{S}^{d-1}$ as follows. 
\rem{}
Fix $0<\gotd<\delta$ and let $C>0, N\in\Bbb{N}$, $c=\{c_n\}_{n\in\Bbb{N}}$ with the same properties as $\tilde{c}=\{\tilde{c}_n\}_{n\in\Bbb{N}}$ in Lemma \ref{Q}.  Let $\delta$ be as in Theorem \ref{main}.
    Take $N_{1}>N\in\Bbb{N}$ large enough such and, for $\b\in\Bbb{S}^{d-1}$  set $\delta_{\beta,0}=\dots=\delta_{\beta,N_1}:=\delta$ and $\bar{\delta}_{\b,0}=\dots=\bar{\delta}_{\b,N_1}=\gotd$, while for $n\in\Bbb{N}$, $n>N_1$  
     define recursively 
   \[
\tilde{\bar{\delta}}_{\b,n+1}:=\max\{\delta\leq \bar{\delta}_{\b,n}:\frac{\Phi(\ell,\b,\delta)}{\|\a\cdot\ell\|} g(c,n,\ell) \leq \max_{\bar{\b}\in\Bbb{S}^{d-1}}\Phi(\ell,\bar{\b},\bar{\delta}_{\bar{\b},n})\quad\forall 0< |\ell|\leq 2^{n+1}\},
\]
with $\Phi$ as in Definition \ref{weakb} 
and $\tilde{\delta}_{\b,n+1}:=\delta_{\b,n}-(\tilde{\bar{\delta}}_{\b,n}-\tilde{\bar{\delta}}_{\b,n+1})$.
From the fact that the fuction $\frac{\max_{\bar{\b}\in\Bbb{S}^{d-1}}\Phi(\ell,\bar{\b},\delta_{\bar{\b}}+\delta')}{\Phi(\ell,\b,\delta+\delta')}$ is an increasing function with respect to $\delta'$, the fact that $\delta>\gotd$ and the fact that $\delta_{\b,n}=\bar{\delta}_{\b,n}+\delta-\gotd$ it follows that
\[
\tilde{\delta}_{\b,n+1}\geq\max\{\delta\leq \delta_{\b,n}:\frac{\Phi(\ell,\b,\delta)}{\|\a\cdot\ell\|} g(c,n,\ell) \leq \max_{\bar{\b}\in\Bbb{S}^{d-1}}\Phi(\ell,\bar{\b},\delta_{\bar{\b},n})\quad\forall 0< |\ell|\leq 2^{n+1}\}.
\]
Then, let:
\[
\delta_{\b,n+1}:=\min_{|\b-\b'|<C_1 2^{dn}e^{-{2^n {c}_n}}}\tilde{\delta}_{n+1,\b'}(1-{c}_{n+1}),
\]
\[
\bar{\delta}_{\b,n+1}:=\min_{|\b-\b'|<C_1 2^{dn}e^{-{2^n {c}_n}}}\tilde{\bar{\delta}}_{n+1,\b'}(1-{c}_{n+1}).
\]
Then, from the fact that $\inf_{\b\in\Bbb{S}^{d-1}}\inf_{n\in\Bbb{N}}\bar{\delta}_{\b,n}>0$ it follows that:
\[
\inf_{n\in\Bbb{N}}\delta_{\b,n}\geq \delta-(\gotd- \inf_{n\in\Bbb{N}}\bar{\delta}_{\b,n})\geq \delta-\gotd.
\]
In particular, for $\gotd$ small enough we have that:
\[
\inf_{n\in\Bbb{N}}\delta_{\b,n}\geq \frac{3\delta}{4}.
\]
\erem

\begin{definition}\label{scoglio}
Fix ¨$\gotd<\delta$ as in the remark above and let $C>0, N\in\Bbb{N}$, $c=\{c_n\}_{n\in\Bbb{N}}$ with the same properties as $\tilde{c}=\{\tilde{c}_n\}_{n\in\Bbb{N}}$ in Lemma \ref{Q}.  Let $\delta$ be as in Theorem \ref{main}.
    Take $N_{1}>N\in\Bbb{N}$ large enough such and, for $\b\in\Bbb{S}^{d-1}$  set $\delta_{\beta,0}=\dots=\delta_{\beta,N_1}:=\delta$, while for $n\in\Bbb{N}$, $n>N_1$  
     define recursively 
   \[
\tilde{\delta}_{\b,n+1}:=\max\{\delta\leq \delta_{\b,n}:\frac{\Phi(\ell,\b,\delta)}{\|\a\cdot\ell\|} g(c,n,\ell) \leq \max_{\bar{\b}\in\Bbb{S}^{d-1}}\Phi(\ell,\bar{\b},\delta_{\bar{\b},n})\quad\forall 2^n< |\ell|\leq 2^{n+1}\},
\]

with $\Phi$ as in Definition \ref{weakb} and
\[
\delta_{\b,n+1}:=\min_{|\b-\b'|<C_1 2^{dn}e^{-{2^n {c}_n}}}\delta_{n+1,\b'}(1-{c}_{n+1}).
\]
Set also $\e_0=\dots=\e_{N_1}=\e$ and, for $n>N_1$ 
  \[
    \delta_{n}:=\min_{\beta\in\Bbb{S}^{d-1}}\delta_{\b,n},\qquad {\e_{n+1}:=e^{-2^n (c_n+c_{n+1})}\e_n+C e^{2^n c_n} 2^{3(n-1)d}\e_{n-1}^2}.
\]
\end{definition}

\rem{gennecoso}
By construction one has
\[
\e_n e^{-2^n c_n}>C 2^{3(n-1)d}\e_{n-1}^2,
\]
Indeed:
\[
C 2^{3(n-1)d}\e_{n-1}^2<e^{2^n c_n}C 2^{3(n-1)d}\e_{n-1}^2<\e_n.
\]

\erem

\lem{staltrolemma}
There exists a choice of $\gotd>0,{c_n}_{n\in\Bbb{N}},C>0,N>0$ such that $\delta_n$ as above. Then:
\[
\inf_{n\in\Bbb{N}}\delta_{n}>\frac{3\delta}{4},
\]
with $\delta$ as in Theorem \ref{main}.
\elem

\proof
It follows directly from the definition of $\{c_n\}_{n\in\Bbb{N}},C,N_1$.

\qed

\lem{questoqui}
Let $\xi$ and $\Bbb{T}^d_\xi$ be as in Definition \ref{param} and consider 
 $f\in D^{\omega}(\Bbb{T}^d_\xi)$.
 For all $N\in\Bbb{N}$, $\b\in\Bbb{S}^{d-1}$ and for all $\bar{\delta}_{\b}<\delta_{\beta}$ such that $|\delta-\bar{\delta}_{\b}|<\frac{\delta}{4}$
 where
  $\delta:=\inf_{\beta\in\Bbb{S}^{d-1}}\delta_{\b}$, one has
\[
|\R_{N}f_{\b}|_{\bar{\delta}_{\beta}}\leq e^{-\frac{(\delta_{\b}-\bar{\delta}_{\b})N}{2}}|\R_{N}f_{\b}|_{\delta_{\beta}}+e^{-\frac{\delta N}{4}}\|\R_N f\|_{\delta}.
\]
\elem

\proof
Recall \eqref{segmento} and \eqref{stanormina}. Then one has
\beq{duesomme}
|\R_{N}f_{\b}|_{\bar{\delta}_{\beta}} \leq  
 \sum_{\substack{\ell\in\Bbb{Z}^d, \, |\ell|>N \\ |\b\cdot\ell|>\frac{N}{2}}}\sup_{z\in S_{\b,\bar{\delta}_{\b}}}\left|\hat{f}(\ell )e^{2\pi i z\cdot\ell}\right|
+\sum_{\substack{\ell\in\Bbb{Z}^d,\,  |\ell|>N \\ |\b\cdot\ell|\leq \frac{N}{2}}}\sup_{z\in S_{\b,\bar{\delta}_{\b}}}\left|\hat{f}(\ell)e^{2\pi i z\cdot\ell}\right|.
\eeq

The first summand in the r.h.s.~of \eqref{duesomme} is bounded by
%5 %%ho aggiunto piu passaggi%%
{
\[
\sum_{\substack{\ell\in\Bbb{Z}^d \\  |\b\cdot\ell|>\frac{N}{2}}}\sup_{x+i\b y\in S_{\b,\bar{\delta}_{\b}}} \left|\hat{f}(\ell)e^{2\pi i (x+i\b y)\cdot\ell}\right|\leq \sum_{\substack{\ell\in\Bbb{Z}^d \\  |\b\cdot\ell|>\frac{N}{2}}} \left|\hat{f}(\ell)e^{2\pi \bar{\delta}_\b|\b\cdot\ell|}\right| 
\leq e^{-\frac{N(\delta_{\b}-\bar{\delta}_{\b})}{2}}|\R_{N}f_{\b}|_{\delta_{\b}},
\]
}
{where the last inequality  follows from the fact that for $\ell\in\Bbb{Z}^d$ such that $ |\b\cdot\ell|>\frac{N}{2}$ we have
$$
\left|\hat{f}(\ell)e^{2\pi \bar{\delta}_\b|\b\cdot\ell|}\right|
\leq e^{-\frac{N(\delta_{\b}-\bar{\delta}_{\b})}{2}} \left|\hat{f}(\ell)e^{2\pi \delta_\b|\b\cdot\ell|}\right|
=e^{-\frac{N(\delta_{\b}-\bar{\delta}_{\b})}{2}} \sup_{x+i\b y\in S_{\b,\delta_{\b}}} \left|\hat{f}(\ell)e^{2\pi i (x+i\b y)\cdot\ell}\right|.
$$
}

{
Regarding the second summand in the r.h.s.~of \eqref{duesomme} we start by noticing that
for $|\ell|>N$ one has 
\begin{equation}\label{stadise}
|\hat{f}(\ell)|\leq e^{-2\pi\delta |\ell|_2}\|\R_{N}f\|_{\delta}.
\end{equation} 
{Indeed using Remark \ref{deltaxi} we may write
$$
\|\R_{N}f\|_{\delta}=\sum_{\substack{\ell\in\Bbb{Z}^d \\ |\ell|>N}}|\hat{f}(\ell)|e^{2\pi|\ell|_2\delta},
$$
so that \eqref{stadise} follows.
}
But then, by Lemma \ref{staltrolemma} 
we bound the second summand in the r.h.s. of \eqref{duesomme} by
\[
\begin{aligned}
\sum_{\substack{\ell\in\Bbb{Z}^d,\,|\ell|>N \\  |\b\cdot\ell|\leq \frac{N}{2}}}
 \sup_{z\in S_{\b,\bar{\delta}_{\b}}} \left|\hat{f}(\ell)e^{2\pi i z\b\cdot\ell}\right| 
&\leq 
\|\R_{N}f\|_{\delta}\sum_{\substack{\ell\in\Bbb{Z}^d,\,|\ell|>N \\  |\b\cdot\ell|\leq \frac{N}{2}}} \left|e^{-2\pi|\ell|_2\delta}e^{2\pi \frac{N}{2}\bar{\delta}_{\b}}\right| \\
&\leq \|\R_{N}f\|_{\delta}\sum_{\substack{\ell\in\Bbb{Z}^d,\,|\ell|>N \\  |\b\cdot\ell|\leq \frac{N}{2}}} \left|e^{-2\pi( |\ell|_2 + \frac{5N}{8})\delta}\right|
\leq e^{-\frac{\delta N}{4}}\|\R_{N}f\|_{\delta},
\end{aligned}
\]
so the assertion follows.}
\qed

The aim of this Section is to prove the following ``KAM step''.

\pro{KAMstep}
Let $\a$ be a weak-Bryuno vector and for all $n\in\Bbb{N}$ set $\xi_n:=\{\delta_{\b,n}\}_{\b\in\Bbb{S}^{d-1}}$ with $\delta_{\b,n}$ as in Definition \ref{scoglio}.
Consider the family of
thickened tori $\Bbb{T}^d_{\xi_n}$, and take any $f_n=R_{\a}+\Delta f_n\in D^{\omega}(\Bbb{T}^d_{\xi_n})$ such that 
\[
\lim_{m\to\infty} \frac{f_n^m}{m}=\a,
\]
satisfying
\beq{stacosa}
\|\mathcal{T}_{2^{n}}\Delta f_{n}\|_{\xi_{n}}<C^{n-1}\e_{n-1}^2,\quad \|\mathcal{R}_{2^{n}}\Delta f_{n}\|_{\xi_{n}}<\e_{n}.
\eeq

Then, there exists $H_{n}\in D^{\omega}(\Bbb{T}^d_{\xi_{n+1}})$ with
\beq{stimah}
\|H_{n}-id.\|_{\xi_{n+1}}<C^{n} \e_n
\eeq
such that
\[
f_{n+1}:=H_{n}\circ f_n\circ H_{n}^{-1}= R_\a + \Delta f_{n+1},
\]
satisfies
\beq{stacosa1}
\|\mathcal{T}_{2^{n+1}}\Delta f_{n+1}\|_{\xi_{n+1}}<C^{n+1}\e_{n}^2,\quad \|\mathcal{R}_{2^{n+1}}\Delta f_{n+1}\|_{\xi_{n+1}}<\e_{n+1}.
\eeq
\epro

\rem{concludo}
Proposition \ref{KAMstep} above implies immediately our main Theorem \ref{main} by setting 
\[
H:=\lim_{N\to\infty} H_N\circ\ldots\circ H_1.
\]
\erem

To prove Proposition \ref{KAMstep} we start by introducing some further notation.

Introduce the projections
\beq{proj}
\begin{aligned}
\pi_i\;:\;&\Bbb{T}^d\longrightarrow [0,1],\\
&x \longmapsto x_i,
\end{aligned}
\qquad i=1,\ldots,d,
\eeq
\beq{proj1}
\begin{aligned}
\tilde{\pi}_i\;:\;&\Bbb{T}^d\longrightarrow \Bbb{R},\\
&x \longmapsto x_i,
\end{aligned}
\qquad i=1,\ldots,d,
\eeq
where we identify $\Bbb{T}^d$ with the fundamental domain $[0,1]^d$.
\lem{AA}
Let $\x=\{\delta_{\b}\}_{\b\in\Bbb{S}^{d-1}}$ be as in Definition \ref{param}, let $f\in{\rm Diff}^{\omega}(\Bbb{T}^{d}_{\xi})$ be such that for every $x\in\Bbb{T}^d$ one has
\[
\lim_{m\to +\infty} \frac{f^m(x)}{m} ({\rm{mod}} \Bbb{Z}^d)=\a ({\rm{mod}} \Bbb{Z}^d),
\]
 and set
 $\Delta f:=f-R_{\a}$. Then, there exists $x^{(1)},\ldots,x^{(d)}\in\Bbb{T}^d$ such that $\pi_j(\Delta f)(x^{(j)})=0$ for $j=1,\ldots,d$.
\elem

\proof
Let $\tilde{f}\in D^{\omega}(\Bbb{T}^d_\x)$ be a lift of $f$ on $\Bbb{R}^d$ such that for $x\in\Bbb{R}^d$, $f(x)=x+\a+\Delta \tilde{f}(x)$ and with $\int_{[0,1]^d}\Delta \tilde{f}(x)dx=0$. 
Then, for all $x\in\Bbb{R}^d$ one has
\[
\lim_{m\to+\infty}\frac{f^m(x)}{m}=\a+\lim_{m\to+\infty} \su{m}\sum_{i=0}^{m-1}\Delta f\circ f^{i}(x),
\]
and hence
\[
\lim_{m\to+\infty} \su{m}\sum_{i=0}^{m-1}\Delta \tilde{f}\circ \tilde{f}^{i}(x)=0.
\]
In particular, for $j=1,\ldots,d$,
\[
\int_{[0,1]^d}\tilde{\pi}_{j}\Delta \tilde{f}(x) dx=\lim_{m\to+\infty} \su{m}\sum_{i=0}^{m-1}\tilde{\pi}_j\Delta \tilde{f}\circ \tilde{f}^{i}(x)=0,
\]
so that for $j=1,\ldots,d$, there exists $x^{(j)}\in\Bbb{T}^d$ such that $\pi_{j}\Delta f(x^{(j)})=0$.
\qed

\lem{resto}
Let $\xi=\{\delta_{\b}\}_{\b\in\Bbb{S}^{d-1}}$ be as in Definition \ref{param},  let $f\in D^{\omega}(\Bbb{T}^{d}_{\xi})$ such that 
\[
\lim_{m\to+\infty} \frac{f^m}{m}=\a,
\]
and set $\Delta f:=f-R_{\a}$. Let $\bar{\x}=\{\bar{\delta}_{\b}\}_{\b\in\Bbb{S}^{d-1}}$ be another family as in Definition \ref{param} such that  $0<\bar{\delta}_{\b}\leq \delta_{\b}$
for all $\b\in\Bbb{S}^{d-1}$. 
Fix $N\in \Bbb{N}$
and suppose that there exists $H_+=id.+h_+\in D^{\omega}(\Bbb{T}^d_{\bar{\xi}})$ with $H_+^{-1}\in D^{\omega}(\Bbb{T}^d_{\bar{\xi}})$ such that
\beq{C}
h_+\circ R_{\a}-h_+=\T_{N}f-\hat{f}(0).
\eeq
Let $g(x):=H_+^{-1}(x)-x+h_+(x)$ and ${f}_+(x)=H_+^{-1}\circ f\circ H_+(x)$.
Then, for all $x\in \Bbb{T}^d_{\bar{\xi}}$ one has
\begin{subequations}\label{perfetto}
\begin{align}
{f}_+(x)&=R_{\a}(x)+\hat{f}(0)+\R_{N}\Delta f(x)+(\Delta f(x+h_+(x))-\Delta f(x)) \\
&\qquad+(h_+(x+\a)-h_+(x+\a+h_+(x)+\Delta f(x+h_+(x)))) \\
&\qquad
+g(x+h_+(x)+\a+\Delta f(x+h_+(x))).
\end{align}
\end{subequations}
\elem

\proof
By direct calculation one finds
\beq{conto}
\begin{aligned}
{f}_+(x)&=H_+^{-1}\circ f\circ H_+(x)=H_+^{-1}(x+h_+(x)+\a+\Delta f(x+h_+(x))) \\
&=x+h_+(x)+\a+\Delta f(x+h_+(x))-h(x+h_+(x) +\a+\Delta f(x+h_+(x)))\\
&\qquad+g(x+h_+(x)+\a+\Delta f(x+h_+(x))) \\
& =x+\a+(h_+(x)-h_+(x+\a)+\Delta f(x))+(\Delta f(x+h_+(x))-\Delta f(x)) \\
&\qquad
+(h_+(x+\a)-h_+(x+\a+h_+(x)+\Delta f(x+h_+(x)))) \\
&\qquad
+g(x+h_+(x)+\a+\Delta f(x+h_+(x))).
\end{aligned}
\eeq
In particular, by \eqref{C}
\[
h_+(x)-h_+(x+\a)+\Delta f(x)=\hat{f}(0)+\R_{N}\Delta f(x),
\]
so that the assertion follows.
\qed

\lem{average}
Under the same hypotheses of Lemma \ref{resto} there exist $x^{(1)},\ldots,x^{(d)}\in\Bbb{R}^d$ such that
\[
\begin{aligned}
\hat{f}_i(0)&=\pi_i(-\R_{N}\Delta f(x^{(i)})-(\Delta f(x^{(i)}+h_+(x^{(i)}))-\Delta f(x^{(i)}))) \\
&\qquad
-\pi_i(h_+(x^{(i)}+\a)-h_+(x^{(i)}+\a+h_+(x^{(i)})+\Delta f(x^{(i)}+h_+(x^{(i)})))) \\
&\qquad+
\pi_i(g(x^{(i)}+h_+(x^{(i)})+\a+\Delta f(x^{(i)}+h_+(x^{(i)})))),
\end{aligned}
\]
for $i=1,\ldots,d$.
\elem

\proof
Let ${f}_+(x):=H_+^{-1}\circ f\circ H_+(x)$. Then
\[
\lim_{m\to+\infty} \frac{{f}_+^m(x)}{m}=\lim_{m\to+\infty}  \frac{f^m(x)}{m}=\a.
\]
By Lemma \ref{AA} there exist $x^{(1)},\ldots,x^{(d)}\in\Bbb{R}^d$ such that for $i=1,\ldots,d$, 
\[
\pi_i({f}_+(x^{(i)})-R_{\a}(x^{(i)}))=\pi_i(\Delta {f}_+(x^{(i)}))=0,
\] 
so the assertion follows by Lemma \ref{resto}.
\qed

%%
%%
%%\lem{sauna}
%%Fix $\b\in[0,1]$ and for all $l\in \Bbb{N}$ let $\delta_{\b,l}$ be as in Definition \ref{scoglio}. Let also $N(l)$ and 
%%$\delta_{l,\beta}^{(i)}$ be as in Definition \ref{musichetta} for $i=0,\ldots,N(l)$.
%%Then, for $|\b-\theta_l|<\su{4\sqrt{l}}$ $\delta_{\b,l}^{(i)}\geq \delta_{\b,l+1}-\frac{16c}{t_n\log t_n}$ and for $|\b-\theta_l|\geq\su{4\sqrt{l}}$ $\delta_{\b,l}^{(i)}\geq \delta_{\b,l+1}-\frac{16c}{t_n^2}$.
%%\elem
%%\proof
%%If $|\b-\theta_l|<\su{4\sqrt{l}}$ we have
%%\[
%%\delta_{\b,l}^{(i)}\geq \delta_{\b,l}-2\sum_{j=0}^{\infty}\frac{4c}{2^jl\log l}=\delta_{\b,l}-\frac{16c}{t_n\log t_n}=\delta_{\b,l+1}+\frac{16c}{t_n\log t_n},
%%\]
%%whereas if $|\b-\theta_l|\geq\su{4\sqrt{l}}$ we have
%%\[
%%\delta_{\b,l}^{(i)}\geq \delta_{\b,l}-2\sum_{j=0}^{\infty}\frac{4c}{2^jl^2}=\delta_{\b,l}-\frac{16c}{t_n^2}=\delta_{\b,l+1}+\frac{16c}{t_n^2}.
%%\]
%%In both cases the assertion follows.
%%\qed

In order to prove Proposition \ref{KAMstep}, we define $H_{n+1}=id.+h_{n+1}$, with $h_{n+1}$ such that
\[
h_{n+1}\circ R_{\a}-h_{n+1}=\T_{2^{n+1}}(f_{n})-\hat{f}_{n}(0),
\]
i.e.~such that
\beq{cena}
h_{n+1}(z):=\sum_{\substack{\ell\in\Bbb{Z}^d\setminus\{0\} \\ |\ell|\leq 2^{n+1} }}\frac{\hat{f}_{n}(\ell)}{e^{2\pi i\a\cdot\ell}-1} e^{2\pi i z\cdot\ell}.
\eeq
Set also
\[
{\phi}_{n+1}(z,\ell):=\frac{\hat{f}_{n}(\ell)}{e^{2\pi i\a\cdot\ell}-1} e^{2\pi i z\cdot\ell}.
\]

%%%
%%%
%%%Then we split $h_{l,i+1}$ as
%%%\[
%%%h_{l,i+1}=h_{l,i+1}^{\paral}+h_{l,i+1}^{\perp},
%%%\]
%%%with
%%%\[
%%%h_{l,i+1}^{\paral}:=\sum_{\substack{\ell\in\Bbb{Z}^2\setminus\{0\} \\ |\ell|\leq 2^{i+1}\tilde{q}_l \\ \ell \paral \nu(l) }}\frac{\hat{f}_{l,i,j}(\ell)}{e^{2\pi i\a\cdot\ell}-1} e^{2\pi i z\cdot\ell},
%%%\]
%%%where $\nu(l)$ is defined in \eqref{nu}, while
%%%\[
%%%h_{l,i+1}^{\perp}:=h_{l,i+1}-h_{l,i+1}^{\paral}.
%%%\]

\lem{3.7}
For any $\b\in\Bbb{S}^{d-1}$ one has
\[
|{\phi}_{n+1,\b}(\cdot,\ell)|_{\tilde{\delta}_{\b,n+1}}\leq C\e_n,.
\]
\elem

\proof
Note that, from the fact that \eqref{stacosa} holds for every $\b\in\Bbb{S}^{d-1}$,
%\[
%|\mathcal{T}_{2^n}\Delta f_{n,\b}|_{\delta_{\b,n}}<C2^{3(n-1)d}\e_{n-1}^2,
%\qquad |\mathcal{R}_{2^n}\Delta f_{n,\b}|_{\delta_{\b,n}}<\e_n
%\]
then for all $\ell\in\Bbb{Z}^{d}$ with $0<|\ell|\leq 2^{n+1}$ and for all $\b\in\Bbb{S}^{d-1}$ we obtain
\[
|\hat{f}(\ell)\Phi(\ell,\b,\delta_{\b,n})|<\e_n,
\]
with $\Phi$ as in Definition $\ref{weakb}$. But then, by definition of $\tilde{\delta}_{\b,n+1}$, for $2^{n}<|\ell|\leq 2^{n+1}$ one has
\[
|{\phi}_{n+1,\b}(\cdot,\ell)|_{\tilde{\delta}_{\b,n+1}}\leq C|\hat{f}(\ell)|\left|\frac{\Phi(\ell,\b,\tilde{\delta}_{\b,n+1})}{\|\a\cdot \ell\|}\right|\leq C\max_{\b\in\Bbb{S}^{d-1}}|\hat{f}(\ell)\Phi(\ell,\b,\delta_{\b,n})|<C\e_n,
\]
while, for $0<|\ell|\leq 2^n$ one has
\[
|{\phi}_{n+1,\b}(\cdot,\ell)|_{\tilde{\delta}_{\b,n+1}}\leq C|\hat{f}(\ell)|\left|\frac{\Phi(\ell,\b,\tilde{\delta}_{\b,n+1})}{\|\a\cdot \ell\|}\right|\leq C2^{3d(n-1)}\e_{n-1}^2 .
\]
In all cases the assertion follows.
\qed

\cor{3.9}
Let $H_{n}=id. + h_{n}$ with $h_{n}$ as in \eqref{cena} and let $\tilde{\xi}_{n+1}:=(\tilde{\delta}_{\b,n+1})_{\b\in\Bbb{S}^{d-1}}$. Then one has
\[
|h_n|_{\tilde{\xi}_{n+1}}<C^{dn}\e_n
\]
\ecor
\proof
By Lemma \ref{3.7}, for all $\b\in\Bbb{S}^{d-1}$ and  all $0<|\ell|\leq 2^{n+1}$ one has
\[
|\phi_{n+1,\b}(\cdot,\ell)|_{\tilde{\delta}_{\b,n+1}}\leq C\e_n.
\]
Then
\[
|h_{n+1}|_{\tilde{\xi}_{n+1}}\leq \sum_{\ell\in\Bbb{Z}^d\setminus\{0\} \\ |\ell|\leq 2^{n+1}}\max_{\b\in\Bbb{S}^{d-1}}|\phi_{n+1,\b}(\cdot,\ell)|_{\tilde{\delta}_{\b,n+1}}\leq C^{n}\e_n,
\]and hence the assertion follows.
\qed

Let $g_{n+1}:=H_{n+1}^{-1}-id.+h_{n+1}$.

\lem{inv} Let $\xi_{n+1}:=\{\delta_{\b,n+1}\}_{\b\in\Bbb{S}^{d-1}}$. Then
$H_{n+1}^{-1}\in D^{\omega}(\Bbb{T}^2_{\xi_{n+1}})$. Moreover one has
\beq{ciotola}
|g_{n+1}|_{\xi_{n+1}}\leq C^{dn}\e_n^2.
\eeq
\elem

\proof
Let $w=H_{n+1}(z)$. Then
\[
g_{n+1}(w)=\int_{0}^{1}Dh_{n+1}(H_{n+1}^{-1}(w)+sh_{n+1}(H_{n+1}^{-1}(w)))h_{n+1}(H_{n+1}^{-1}(w))ds.
\]
In particular, by Corollary \ref{3.9} and the definition of $\xi_{n+1}$, for $z\in \Bbb{T}^d_{\xi_{n+1}}$ all compositions above are well defined: 
indeed by Corollary \ref{3.9}, the image of $\Bbb{T}^d_{\xi_{n+1}}$ under $H_{n+1}$ is contained in $\Bbb{T}^d_{\bar{\xi}_{n+1}}$, 
where $\bar{\xi}_{n+1}=\{\bar{\delta}_{\b,n+1}\}_{\b\in\Bbb{S}^{d-1}}$ is such that
\[
\bar{\delta}_{\b,n+1}:=\max_{|\b-\b'|<C2^{dn}\e_n}\delta_{\b',n+1}(1-C\e_n).
\]

{
Let us now prove the bound \eqref{ciotola}. Since
\[
(id.-h_{n+1}+g_{n+1})\circ(id.+h_{n+1})=id.,
\]
then for any $z\in \Bbb{T}^d_{\xi_{n+1}}$ one has
\[
g_{n+1}(z+h_{n+1}(z))=h_{n+1}(z+h_{n+1}(z))-h_{n+1}(z).
\]
Let $w_0(z):=h_{n+1}(z+h_{n+1}(z))-h_{n+1}(z)$ and for $k\geq 1$ let
\[
w_{k}(z):=w_{k-1}(z+h_{n+1}(z))-w_{k-1}(z).
\]
Then one has
\begin{equation}\label{pastalforno}
g_{n+1}(z+h_{n+1}(z))=\sum_{k\geq 0}w_k(z+h_{n+1}(z)),
\end{equation}
provided that the sum converges. Indeed, by definition we have
\[
\begin{aligned}
w_{0}(z+h_{n+1}(z))&=(w_0(z+h_{n+1}(z))-w_0(z))+w_0(z)\\
&=(w_0(z+h_{n+1}(z))-w_0(z))+h_{n+1}(z+h_{n+1}(z))-h_{n+1}(z),
\end{aligned}
\]
%so that $w_0$ is of the solution of the equation that we are considering up to the error $w_0(z+h_{n+1}(z))-w_0(z)$. Then, suppose by induction that:
Now if we assume inductively
\[
\sum_{k=0}^{l}w_{k}(z+h_{n+1}(z))=h_{n+1}(z+h_{n+1}(z))-h_{n+1}(z)+w_{l}(z+h_{n+1}(z))-w_{l}(z),
\]
then
\[
\begin{aligned}
\sum_{k=0}^{l+1}w_{k}(z+h_{n+1}(z))&=h_{n+1}(z+h_{n+1}(z))-h_{n+1}(z)+w_{l}(z+h_{n+1}(z))\\
&\qquad-w_{l}(z)+w_{l+1}(z+h_{n+1}(z)) \\
&=h_{n+1}(z+h_{n+1}(z))-h_{n+1}(z)+w_{l+1}(z+h_{n+1}(z))-w_{l+1}(z),
\end{aligned}
\]
with the last inequality following from the definition of $w_{l+1}(z)$.
In particular, we obtain \eqref{pastalforno}
provided that the sum converges.

Now, for $k\geq 0$, let $w_{k,1},\dots,w_{k,d},h_{n+1,1},\ldots,h_{n+1,d}\in C^{\omega}_{\xi_{n+1}}(\Bbb{T}^d)$ be such that $w_{k}(z)=(w_{k,1}(z),\dots,w_{k,d}(z))$
and $h_{n+1}(z)=(h_{n+1,1}(z),\dots,h_{n+1,d}(z))$. Then, for $j=1,\dots, d$ one has
\[
w_{k+1,j}(z)=w_{k,j}(z+h_{n+1}(z))-w_{k,j}(z)=\sum_{t_1,\dots,t_d\in\Bbb{N}}\frac{\partial^{t_1+\dots+t_d}w_{k,j}}{\partial_1^{t_1}\ldots\partial_d^{t_d}}\frac{h_{n+1,1}^{t_1}\ldots h_{n+1,d}^{t_d}}{t_1!\ldots t_d!}
\]
so that iterating we obtain
\[
w_{k+1,j}(z)=\sum_{\substack{t_{1},\dots,t_{d}\in\Bbb{N} \\ t_1+\ldots + t_d \geq k+1}}\frac{\partial^{t_1+\ldots+t_{d}}h_{n+1,j}}{\partial_1^{t_1}\ldots\partial_d^{t_{d}}}\frac{h_{n+1,1}^{t_1}\dots h_{n+1,d}^{t_{d}}}{t_1!\ldots t_{d}!}
\]

By Lemma \ref{prodotto} and Lemma \ref{derivate} we have
\[
\begin{aligned}
\|w_{k+1}\|_{\xi_{n+1}}&\leq \sum_{\substack{t_1,\dots,t_{d}\in\Bbb{N} \\ t_1+\ldots+t_d\geq k+1}}\frac{(4\pi)^{(n+1)d(t_1+\dots+t_{d})} 
\|h_{n+1}\|_{\xi_{n+1}}^{t_1+\dots t_{d}}}{t_1!\dots t_{d}!} \\
&=\sum_{m\geq k+1}\sum_{t_1+\dots +t_{d}=m}\frac{(4\pi)^{(n+1)dm}\|h_{n+1}\|_{\xi_{n+1}}^{m+1}}{t_1!\dots t_{d!}}\leq C^{(k+1)d}\|h_{n+1}\|_{\xi_{n+1}}^{k+2},
\end{aligned}
\]
which in turn implies
\[
\|g_{n+1}\|_{\xi_{n+1}}\leq \sum_{k\geq 0} \|w_k\|_{\xi_{n+1}}\leq C^{nd}\|h_{n+1}\|_{\xi_{n+1}}^2.
\]
}

In particular, frome the definition of $\xi_{n+1}$ and the fact that $\e_n=Ce^{-2^{n}c_n}\e_{n-1}$ we deduce that $\Bbb{T}^d_{\bar{\xi}_{n+1}}\subset \Bbb{T}^d_{\tilde{\xi}_{n+1}}$,
and hence, $H^{-1}_{n+1}\in D^{\omega}(\Bbb{T}^d_{\xi_{n+1}})$. Moreover one has
\[
|Dh_{n+1}|_{\tilde{\delta}_{\b,n+1}}\leq C 2^{2nd}\e_n
\]
so that, using also Corollary \ref{3.9}, we obtain
\beq{3.51}
\begin{aligned}
\|g_{n+1}\|_{\xi_{n+1}}&\leq |Dh_{n+1}|_{\tilde{\delta}_{\tilde{\xi}_{n+1}}}\|h_{n+1}\|_{\tilde{\xi}_{n+1}} \\
&\leq C2^{3dn}\e_{n}^2
\end{aligned},
\eeq
from which \eqref{ciotola} follows.
\qed

\lem{tantoper}
For all $n\in\Bbb{N}$ one has
\[
|\hat{f}_{n+1}(0)|\leq C2^{2dn}\e_{n}^{2}.
\]
\elem

\proof
It follows directly by Lemma (\ref{average}) combined with Corollary \ref{3.9}.
 \qed

\giu
\nl
{\bf{Proof (of Proposition \ref{KAMstep})}} 
Thanks to Lemmata above we are only left to prove the bounds in \eqref{stacosa1}.
%
%\beq{}
%\begin{aligned}
%&|T_{2^{i+1}\tilde{q}_l}\Delta f_{l,i+1,\beta}-\Delta \hat{f}_{l+1,i}(0)|_{\delta_{\b,l}^{(i+1)}}\leq \e_{\b,l,i+1}, \\
%&|R_{2^{i+1}\tilde{q}_l}\Delta f_{l,i+1,\b}|_{\delta_{\b,l}^{(i+1)}}\leq \e_{\b,l,i+1}.
%\end{aligned}
%\eeq
%
%
%\beq{}
%|T_{2^{i+1}\tilde{q}_l}\Delta f_{\a,l,i+1}-\Delta \hat{f}_{\a,l,i+1}(0)|_{\delta_{\a,l}^{(i+1)}}\leq 
%\eeq
%\beq{}
%\sup_{|\a-\b|<\su{2\sqrt{t_n}(1+\su{2^i})}}|\Delta f_{\a,l,i}|_{\delta_{\a,l}^{(i)}}|h_{\b,l,i}|_{\delta_{\b,l}^{(i+1)}}+\sup_{|\a-\b|<\su{2\sqrt{t_n}(1+\su{2^i})}}|\Delta f_{\b,l,i}|_{\delta_{\b,l}^{(i)}}|h_{\a,l,i}|_{\delta_{\a,l}^{(i+1)}}\leq \e_{\a,l,i+1}
%\eeq
%and
%\beq{}
%|R_{2^{i+1}\tilde{q}_l}\Delta f_{\a,l,i+1}|_{\delta_{\a,l}^{(i+1)}}\leq e^{-(\delta_{\a,l}^{(i)}-\delta_{\a,l}^{(i+1)})2^i\tilde{q}_l}|R_{2^{i+1}\tilde{q}_l}\Delta f_{\a,l,i+1}|_{\delta_{\a,l}^{(i)}}<\e_{\a,l,i+1}.
%\eeq
%\elem
By Lemma \ref{inv} we see that $f_{n+1}:=H_{n+1}^{-1}\circ f_{n}\circ H_{n+1}\in D^{\omega}(\Bbb{T}^2_{\x_{n+1}})$ and by Lemma \ref{resto} we know that

\begin{subequations}\label{nonloso}
\begin{align}
\Delta f_{n+1}(x)-&\Delta\hat{f}_{n+1}(0)=\R_{2^{n+1}}\Delta f_{n}(x)+(\Delta f_{n}(x+h_{n+1}(x))-\Delta f_{n}(x)) \\
&+(h_{n+1}(x+\a)-h_{n+1}(x+\a+h_{n+1}(x)+\Delta f_{n}(x+h_{n+1}(x)))) \\
&+g_{n+1}(x+h_{n+1}(x)+\a+\Delta f_{n}(x+h_{n+1}(x))),
\end{align}
\end{subequations}
see \eqref{perfetto}. Now, by Lemma \ref{questoqui} and the fact that for each $l\in\Bbb{N}$ one has $\delta_{\b,l}\geq\frac{\delta}{2}$, we deduce
\[
\begin{aligned}
|\R_{2^{n+1}}\Delta f_{n}|_{\delta_{\b,n+1}} &\leq e^{-\frac{(\delta_{\b,n}-\delta_{\b,n+1})2^{n+1}}{2}}|\R_{2^{n+1}}\Delta f_{n}|_{\delta_{\b,n+1}}\\
&\qquad\qquad+e^{-\frac{\delta 2^{n+1}}{2}} \|\R_{2^{n+1}}\Delta f_{n}\|_{\frac{\delta}{2}} \\
&\leq \e_n^2+e^{-\delta 2^n}\e_n  ,
\end{aligned}
\]
On the other hand by Lemma \ref{3.7} we deduce

\beq{suppli}
\begin{aligned}
|(\Delta f_{n,\b}(x+h_{n+1}(x))-&\Delta f_{n,\b}(x))|_{\delta_{\b,n+1}} = \\
&=
\left|  \int_{0}^{1}D\Delta f_{n,\b}(x+sh_{n+1}(x))h_{n+1,\b}(x)ds \right|_{\delta_{\b,n+1}} \\
&\leq \Big(\sup_{ \substack{\b'\in\Bbb{S}^{d-1}\\ |\b-\b'|<C2^{dn}\e_n}}|D f_{n,\b'}|_{\delta_{n+1,\b'}}\Big)|h_{n+1,\b}|_{\delta_{\b,n+1}} \\
&<C2^{3dn}\e_n.
\end{aligned}
\eeq

Similarly
\beq{doppiosuppli}
\begin{aligned}
|h_{n+1,\b}&(x+\a)-h_{n+1,\b}(x+\a+h(x)+\Delta f_{n}(x+h_{n+1}(x)))|_{\delta_{\b,n+1}} \\
&\leq  \int_{0}^{1} \Big| Dh_{n+1,\b}(x+\a+s(h_{n}(x)+ \\ 
&\qquad\quad+\Delta f_{n}(x+h_{n+1}(x))))(h_{n+1,\b}(x)+\Delta f_{n,\b}(x+h_{n}(x)))  \Big|_{\delta_{\b,n+1}} ds \\
&\leq \Big(\sup_{\substack{\b'\in\Bbb{S}^{d-1} \\ |\b-\b'|<C2^{dn}\e_n}}  |D h_{n+1,\b'}|_{\delta_{n+1,\b'}}\Big) ( |h_{n+1,\b}|_{\delta_{\b,n+1}}  + |\Delta f_{n,\b}|_{\delta_{\b,n+1}} )\\
&<C2^{3dn}\e_n.
\end{aligned}
\eeq
Finally, by Lemma \ref{inv} we have
\beq{crocchetta}\begin{aligned}
|g_{n+1,\b}(x+h_{n+1}(x)+\a+&\Delta f_{n}(x+h_{n+1}(x)))|_{\delta_{\b,n+1}} \leq \\
&\leq \sup_{\substack{\b'\in\Bbb{S}^{d-1} \\ |\b-\b'|<C2^{dn}\e_n}}  |g_{n+1,\b'}|_{\delta_{n+1,\b'}} \\
&<C2^{3dn}\e_n^2.
\end{aligned}
\eeq

Combining \eqref{suppli}, \eqref{doppiosuppli} and \eqref{crocchetta} we obtain
\[
|f_{n+1}|_{\xi_{n+1}}<C2^{3nd}\e_n^2<\e_{n+1}.
\]
In particular, \eqref{stacosa} follows.
\qed

%%%%%%%%%%%%%%%%%%%%%%%%%%%%%%%%%%%%%%%%%%%%%%%%%%%%%%%%%%%%%%%%%%%%%%%%% 
\appendix 
%%%%%%%%%%%%%%%%%%%%%%%%%%%%%%%%%%%%%%%%%%%%%%%%%%%%%%%%%%%%%%%%%%%%%%%%% 
 
%%%%%%%%%%%%%%%%%%%%%%%%%%%%%%%%%%%%%%%%%%%%%%%%%%%%%%%%%%%%%%%%%%%%%%%%% 
%%%%%%%%%%%%%%%%%%%%%%%%%%%%%%%%%%%%%%%%%%%%%%%%%%%%%%%%%%%%%%%%%%%%%%%%% 
\zerarcounters 
\section{On the one dimensional case} 
\label{appendice} 
%%%%%%%%%%%%%%%%%%%%%%%%%%%%%%%%%%%%%%%%%%%%%%%%%%%%%%%%%%%%%%%%%%%%%%%%% 
%%%%%%%%%%%%%%%%%%%%%%%%%%%%%%%%%%%%%%%%%%%%%%%%%%%%%%%%%%%%%%%%%%%%%%%%% 

Here we prove that a number is a weak-Bryuno number if and only if it is a Bryuno number. 
 %In order to make the definition more clear now we will show that in dimension one the two definitions are equivalent. Indeed, l

Assume first $\a$ to be a Bryuno number. Fix $\gotd>0$ $C=1$ and, for $n\in\Bbb{N}$ set 
\[
c_n:=2\pi\frac{1}{2^n} \log\su{\Omega_{\a}(2^n)}
\]
 so that $\{c_n\}_{n\in\Bbb{N}}$ is in $\ell^1(\Bbb{R}^{+})$. Fix also $N\in\Bbb{N}$ such that:
 \[
 \gotd>\sum_{n\geq N}\frac{1}{2^n}\log\su{\Omega_{\a}(2^n)}.
 \] 
 
 Note that for $d=1$ in Definition \ref{weakb}, then $\Bbb{S}^0:=\{1,-1\}$ and for any $\ell\in\Bbb{Z}$, $\delta>0$ one has
 \[
 \Phi(\ell,1,\delta)=\Phi(\ell,-1,\delta)=e^{2\pi|\ell|\delta},
 \]
 which in turn means that for any $n\ge 0$ one has
 $\delta_{1,n}=\tilde{\delta}_{1,n}=\tilde{\delta}_{-1,n}$ and hence $\delta_{1,n}=\delta_{-1,n}=:\delta_n$.
 
 %Let $n\geq 0$ and suppose that $\delta_n$ is defined as in Definition \ref{weakb} (note that if $d=1$, then $\delta_{n,1}=\tilde{\delta}_{1,n}=\tilde{\delta}_{-1,n}$ and $\delta_{1,n}=\delta_{-1,n}$). 

Now, for any $\delta>0$, if $\ell\in\Bbb{Z}$ is such that $2^{n}<|\ell|\leq 2^{n+1}$, we have
$$
\frac{\Phi(\ell,1,\delta)}{\|\a\ell\|}g(c,n,\ell)=\frac{\Phi(\ell,1,\delta)}{\|\a\ell\|}\leq 
e^{2\pi(\delta+\su{2\pi2^n}\log\su{\Omega_{\a}(2^{n+1})})|\ell|},
$$ 
while, if $0<|\ell|\leq 2^{n}$ we have
$$
\frac{\Phi(\ell,1,\delta)}{\|\a\ell\|}g(c,n,\ell)
\leq e^{2\pi|\ell|\delta}e^{(1-2\pi)\log\su{\Omega_{\a}(2^n)}}
\leq e^{2\pi\delta|\ell|}.
$$

%%$$
%%\frac{\Phi(\ell,1,\delta)}{\|\a\cdot\ell\|}g(c,n,l)\leq e^{2\pi|l|(\delta+\frac{\log\Omega_{\a}{\su{|l|}}}{|l|}-\frac{2^n\log\su{\Omega_{\a}(2^n)}}{|l|})}\leq e^{2\pi\delta|l|}
%%$$
This means that
\[
\delta_{n+1}\geq\delta_n-\su{ 2^{n+1}}\log\su{\Omega_{\a}(2^{n+1})},
\]
which in turn implies
\[
\inf_{n\in\Bbb{N}}\delta_n\ge\gotd-\sum_{n\geq N}\frac{1}{2^n}\log\su{\Omega_{\a}(2^n)}>0,
\]
i.e.~$\a$ is a weak-Bryuno number. 

On the other hand by definition if $\a$ is not a Bryuno number, then it is not a weak-Bryuno number. However a stronger result holds. In fact, if 
$\a$ is not a Bryuno number,
%Now let us show that if $\a$ is not a Bryuno number, then it is not a weak-Bryuno number. 
let us fix any $C>0$, $N\in\Bbb{N}$ and $c\in \ell^1(\Bbb{R}^{+})$. Let $n_k$ be an increasing sequence of positive integers with $n_0>N$ and $2^{n_k}<l_k\leq 2^{n_{k+1}}$ such that
\begin{itemize}
    \item For $k\in\Bbb{N}$ we have $\Omega_{\a}(2^{n_k+1})=\|l_k\a\|$,
    \item  $\sum_{k\in\Bbb{N}}\frac{1}{2^{n_k}}\log\su{\Omega_{\a}(2^{n_k})}=\infty$.
\end{itemize}
Then for any $\delta>0$ one has
$$
\frac{\Phi(l_k,1,\delta)}{\|l_k\a\|}\geq e^{2\pi|l_k|(\delta+\frac{1}{2^{n_k+1}}\log\su{\Omega_{\a}(2^{n_k})})},
$$
so in particular one has 
$$
\delta_{n_k+1}\leq \delta_{n_k}-\frac{1}{2^{n_k+1}}\log\su{\Omega_{\a}(2^{n_k})},
$$
which in turn implies 
$$
\inf_{n\in\Bbb{N}}\delta_{n}\leq \delta_0-\sum_{k\in\Bbb{N}}\frac{1}{2^{n_k+1}}\log\su{\Omega_{\a}(2^{n_k})}=-\infty.
$$ 
Thus if $\a$ is not a Bryuno number, then \eqref{weakcond} is not satisfied as well.
%In particular, if $\a$ is not a Bryuno number, then it is not a weak-Bryuno number, from which we get that the two definitions are equivalent for $d=1$. 

%%%%%%%%%%%%%%%%%%%%%%%%%%%%%%%%%%%%%%%%%%%%%%%%%%%%%%%%%%%%%%%%%%%%%%%%%%
%%%%%%%%%%%%%%%%%%%%%%%%%%%%%%%%%%%%%%%%%%%%%%%%%%%%%%%%%%%%%%%%%%%%%%%%%%
 

\begin{thebibliography}{99} 
%%%%%%%%%%%%%%%%%%%%%%%%%%%%%%%%%%%%%%%%%%%%%%%%%%%%%%%%%%%%%%%%%%%%%%%%%%
%%%%%%%%%%%%%%%%%%%%%%%%%%%%%%%%%%%%%%%%%%%%%%%%%%%%%%%%%%%%%%%%%%%%%%%%%%
{\small 


\bibitem{AK}{
D.~Anosov, A.~Katok, {\sl New examples in smooth ergodic theory, ergodic diffeomorphisms.} 
Trans. Mosc. Math. Soc, 23(1), 1970.
}

\bibitem{avilaglobal}{A.~Avila, {\sl Almost reducibility and absolute continuity I},  preprint 2010 arXiv:1006.0704}

\bibitem{avila_global}{ 
A.~Avila, 
{\sl KAM, Lyapunov exponents, and the Spectral Dichotomy for typical one-frequency Schrodinger operators.} 
 preprint 2023, arXiv:2307.11071.
}

\bibitem{BF}{
A.~Bounemura, J.~F\'ejoz,
\textit{KAM, $\alpha$-Gevrey regularity and the $\alpha$-Bruno-R\"ussmann condition}
Ann. Sc. Norm. Super. Pisa Cl. Sci. (5) {\bf{19}} (2019), no.~4, 1225--1279.
}

\bibitem{Bounemoura}{
A.~Bounemoura,
\textit{Some remarks on the optimality of the Bruno-R\"ussmann condition},
Bull. Soc. Math. France \textbf{147} (2019), no.~2, 341--353.
}

\bibitem{Bounemoura 1}{
A.~Bounemoura,
\textit{Optimal linearization of vector fields on the torus in non-analytic Gevrey classes}. 
Ann. Inst. H. Poincar\'e C Anal. Non-Lin\'eaire 39(2022), no.3, 501--528. %501-528.
}




\bibitem{Bry}{
A.D.~Bryuno,
\textit{Analytic form of differential equations. I, II},
Trudy Moskov. Mat. Ob\v{s}\v{c}.
\textbf{25} (1971), 119--262;
ibid. \textbf{26} (1972), 199--239.
English translations:
Trans. Moscow Math. Soc.
\textbf{25} (1971), 131--288 (1973);
ibid. 26 (1972), 199--239 (1974). 
}

\bibitem{BC}{
X.~Buff, A~Ch\'eritat, 
\textit{Upper bound for the size of quadratic Siegel disks.} Invent. Math.
156(1), 1–24 (2004)}

\bibitem{BC1}{ 
X.~Buff, A~Ch\'eritat, 
\textit{The Brjuno function continuously estimates the size of quadratic Siegel disks.} 
Ann. Math. (2) 164(1), 265–312 (2006)
}


\bibitem{C}{
T.~Carletti, {\sl The 1/2-complex Bruno function and the Yoccoz function: a numerical study
of the Marmi–Moussa–Yoccoz conjecture.} Exp. Math. 12(4), 491–506 (2003)
}

\bibitem{Ch}{
A.~Ch\'eritat, {\sl Sur l’implosion parabolique, la taille des disques de Siegel et une conjecture de
Marmi, Moussa et Yoccoz.} Habilitation \`a diriger des recherches, Universit\'e Paul Sabatier
(2008)
}


\bibitem{CC}{
D.~Cheraghi, A.~Ch\'eritat, 
\textit{A proof of the Marmi-Moussa-Yoccoz conjecture for rotation numbers of high type}
Invent. Math. \textbf{202} (2015), no. 2, 677--742.
}


\bibitem{CGP}{
L.~Corsi, G.~Gentile, M.~Procesi, 
\textit{Almost-periodic solutions to the NLS equation with smooth convolution potentials}
preprint 2023.
}

\bibitem{FK}{
B.~Fayad, A.~Katok. {\sl Analytic uniquely ergodic volume preserving maps on odd
spheres.} 
Commentarii Mathematici Helvetici, 89(4):963–977, 2014.
}

\bibitem{GU}{
S.~Ghazouani, C.~Ulcigrai, {\sl A priori bounds for GIETs, affine shadows and rigidity of foliations in genus 2}, 
preprint 2023 arXiv:2106.03529.
}

\bibitem{GU1}{
S.~Ghazouani, C.~Ulcigrai, {\sl Regularity of conjugacies of linearizable generalized interval exchange transformations}, 
preprint 2023 arXiv:2304.06628. 
}

\bibitem{H2}{
M. Herman, {\sl Sur la conjugaison differentiable des diffeomorphismes du cercle a des rotations.} 
Inst. Hautes Etudes Sci. Publ. Math. No. 49 (1979), 5-233.
}


\bibitem{H}{
 M.~Herman, \textit{Une methode pour minorer les exposants de Lyapounov et quelques exemples
 montrant le caractere local d'un theoreme d'Arnold et de Moser sur le tore de dimension 2.}
Comment. Math. Helv. \textbf{58} (1983), no. 3, 453--502.
}


\bibitem{M}{
S.~Marmi, {\sl Diophantine conditions and estimates for the Siegel radius: analytical and
numerical results.}, Nonlinear Dynamics (Bologna, 1988). World Scientific Publishing,
Teaneck, pp. 303–312 (1989)
}


\bibitem{MMY}{
S.~Marmi, P.~Moussa, J.-C.~Yoccoz,\textit{
The Brjuno functions and their regularity properties.} Comm. Math. Phys. \textbf{186} (1997), no.2, 265--293.
}

\bibitem{MMY1}
{S.~Marmi, P.~Moussa, J.-C.~Yoccoz,
\textit{
Complex Brjuno functions.} J. Amer. Math. Soc. \textbf{14} (2001), no.4, 783--841. }




\bibitem{MJ}{
R.~Johnson, J.~Moser, \textit{The rotation number for almost periodic potentials.} Comm. Math.
Phys., \textbf{84} (1982), no. 3, 403--438. 
}

\bibitem{kol}
A.N.~Kolmogorov, 
\emph{On the Conservation of Conditionally Periodic Motions under Small Perturbation of the Hamiltonian}, Dokl. akad. nauk SSSR, 
1954, vol. 98, pp. 527–530. Engl. transl.: 
\emph{Stochastic Behavior in Classical and Quantum Hamiltonian Systems}, Volta Memorial conference, Como, 1977, Lecture Notes in Physics, vol. 93, Springer, 1979, pp. 51--56.


\bibitem{KQY} R.~Krikorian, J.~Wang, J.~You, Q.~Zhou, {\sl Linearization of quasiperiodically forced circle flows beyond Brjuno condition},
 Commun. Math. Phys. 358(1), 81–100 (2018).

\bibitem{Rus}{
H.~R\"ussmann, 
\textit{Nondegeneracy in the perturbation theory of integrable dynamical systems}, 
Stochastics, algebra and analysis in classical and quantum dynamics (Marseille, 1988), 
Math. Appl. \textbf{59},  211--223, Kluwer Acad. Publ., Dordrecht, 1990.
}

\bibitem{Rus1}{H.~R\"ussmann,
\textit{ Stability of elliptic fixed points of analytic area-preserving mappings under the Bruno condition},
Ergodic Theory Dynam. Systems \textbf{22} (2002), no. 5, 1551--1573.
}

\bibitem{sie}{
C.L.~Siegel, {\sl Iteration of analytic functions}, Ann. Math., 1 {\bf 28}, 607 (1942)
}

\bibitem{Y1}{
J.-C.~Yoccoz,
\textit{Lin\'earisation des germes de diff\'eomorphismes holomorphes de $(\Bbb{C},0)$}. 
C. R. Acad. Sci. Paris S\'er. I Math. \textbf{306} (1988), no.~1, 55--58.
}

\bibitem{Y2}{
J.-C.~Yoccoz, 
\textit{Th\'eor\`eme de Siegel, nombres de Bruno et polyn\^omes quadratiques.}
Petits diviseurs en dimension 1.
Ast\'erisque No. \textbf{231} (1995), 3--88.
}


\bibitem{Y3}{
J.-C.~Yoccoz, \textit{Analytic linearization of circle diffeomorphisms.} Dynamical systems
and small divisors (Cetraro, 1998), 125-173, Lecture Notes in Math., 1784, Fond.
CIME/CIME Found. Subser., Springer, Berlin, 2002.
}	


\bibitem{U}{
C.~Ulcigrai. 
{\sl Dynamics and ’arithmetics’ of higher genus surface flows.} 
Proceedings of the International Congress
of Mathematics in Saint Petersburg, 2022.
}


}
\end{thebibliography}
\end{document}